\documentclass[10pt]{article}

\usepackage{amsmath,latexsym,amssymb,graphicx,a4}

\newcommand{\N}        {\mathbb N}
\newcommand{\R}        {\mathbb R}

\newcommand{\Z}        {\mathbb Z}
\newcommand{\E}        {\mathbb E}
\newcommand{\AAA}          {\mathcal{A}} 
\newcommand{\CC}          {\mathcal{C}} 
\newcommand{\BB}          {\mathcal{B}} 
\newcommand{\EE}          {\mathcal{E}}
\newcommand{\II}          {\mathcal{I}} 
\newcommand{\GG}          {\mathcal{G}} 
\newcommand{\LL}         {\mathcal{L}}
\newcommand{\eps}    {\varepsilon} 
\newcommand{\eqdef}     {\stackrel{\triangle}{=}} 
 
\newcommand{\Liminf}    {\mathop{\underline{\rm lim}}}
\newcommand{\proof}        {\paragraph{Proof}}
\newcommand{\diam}      {\hbox{{\rm diam}}}
\newcommand{\carre}     {\hfill$\Box$}
\newcommand{\Var}       {\hbox{{\rm Var}}}

\newtheorem{theorem}      {Theorem}[section]
\newtheorem{definition}      [theorem]{Definition}

\renewcommand{\phi}{\varphi}

\title{Isoperimetry and heat kernel decay on percolation clusters}
\author{Pierre MATHIEU\footnote{CMI, 39 rue Joliot-Curie, 13013 Marseille, FRANCE, 
pierre.mathieu@cmi.univ-mrs.fr}, Elisabeth REMY\footnote{IML, campus de Luminy, 13009 Marseille, FRANCE, 
remy@iml.univ-mrs.fr}} 

\date{}
\begin{document}

\maketitle

\noindent
{\bf Short title}: Isoperimetry on percolation clusters. 

\noindent
{\bf Abstract}: we prove that the heat kernel on the infinite 
Bernoulli percolation cluster in $\Z^d$ almost surely decays faster than $t^{-d/2}$. 
We also derive estimates on the mixing time for the 
random walk confined to a finite box. Our approach is based on local 
isoperimetric inequalities. 
Some of the results of this paper were previously announced in the note \cite{mathieu01a}. 



\section{Introduction}

We deal separately with $2D$ site percolation and bond percolation in any dimension. 

\subsection{Site percolation in $2D$}
\label{intro-2d}

Let $\omega$ be the random sub-graph of $\Z^2$ obtained by keeping (resp. deleting) a point with
probability $p$ (resp. $1-p$), independently for different points of $\Z^2$. Call $Q$ the law 
of $\omega$. Points belonging to $\omega$ are called {\it open}. Two neighbouring points 
of $\omega$ form an {\it open} edge. 
Let $\CC$ denote the open cluster of the origin, {\it i.e.} the connected component of
$\omega$ containing $0$, and define
\begin{displaymath}
 p_c = \sup \{p\,;\,Q[\# \CC=+\infty]= 0 \}\,.
 \end{displaymath}
For a given choice of $\omega$, 
we shall consider the usual random walk on $\CC$, say $(X_t,t\geq 0)$: 
the random walker waits for an exponential time of parameter $1$ and
 then chooses, uniformly at random,  
one of its neighbors in $\Z^2$, say $y$. If 
$y$ belongs to $\CC$ (the edge leading to $y$ is open), 
then the walker moves to $y$; otherwise it stays still. 
Thus $X_t$ defines 
a Markov chain on $\CC$ which is reversible with respect to the counting measure on $\CC$. 
\begin{theorem}
\label{th1}
For any $p>p_c$, there exists a constant $c_1=c_1(p)$ such that $Q$ a.s. 
on the set $\# \CC=+\infty$, and for large enough $t$, we have
\begin{equation}
\label{ineq-th1}
 \sup_{y \in \CC} P_0^{\omega}[X_t=y] 
 \leq 
 \frac{c_1}{t}\,.
\end{equation}
\end{theorem}

\subsection{Bond percolation in $\Z^d$}
\label{intro-d}

Let $\omega$ be the random sub-graph of $\Z^d$ obtained by keeping (resp. deleting) an edge with
probability $p$ (resp. $1-p$) independently for each bond. 
More precisely, for $x, y \in \Z^d$, we write: $x \sim y$ if $x$ and $y$ are neighbors in $\Z^d$, and 
$\E_d=\{(x,y) \in \Z^d \times \Z^d\,,\,x \sim y\}$. 
We identify a sub-graph of $\Z^d$ with an  application $\omega:\E_d \rightarrow \{0,1\}$, 
writing $\omega(x,y)=1$ if the edge $(x,y)$ is present in $\omega$ and writing 
$\omega(x,y)=0$ otherwise. Edges in $\E_d$ will be called {\it open}. 
Let $Q$ be the probability measure on $\{0,1\}^{\E_d}$ under which  
the random 
variables $(\omega(e),\,e \in \E_d)$ are Bernouilli$(p)$ independent variables. 
As before, let 
$p_c = \sup \{p\,;\,Q[\# \CC=+\infty]= 0 \}$ be the critical probability. 

For a given sub-graph $\omega$, let $\CC$ denote the open cluster at the 
origin, {\it i.e.} $\CC$ is the connected component of $\omega$ that contains $0$. 
We still use the notation $X_t$ to denote the random walk on $\CC$: 
the random walker waits for an exponential time of parameter $1$ and
 then chooses, uniformly at random,  
one of its neighbors in $\Z^d$, say $y$. If 
$\omega(x,y)=1$, then the walker moves to $y$; otherwise it stays still. 
Thus $X_t$ defines 
a Markov chain on $\CC$ which is reversible with respect to the counting measure on $\CC$. 
Let $P_x^\omega$ be the law of the chain $(X_t,t\geq0)$ when started at point $x$. 

\begin{theorem}
\label{th2}
For any dimension $d \geq 2$,  for any $p >p_c$, 
there exists a constant $c_1=c_1(p,d)$ such that, $Q$ a.s. on the set $\# \CC=+\infty$, and 
for large enough $t$, we have
\begin{equation}
\label{ineq-th2}
 \sup_{y \in \CC} P_0^{\omega}[X_t=y] 
 \leq 
 \frac{c_1}{t^{d/2}}\,.
\end{equation}
\end{theorem}

\paragraph{Remarks:} 
\begin{itemize}
\item[$\bullet$] it is known that $X_t$ satisfies a central limit theorem, 
see \cite{demasiferrari}. Our estimate (\ref{ineq-th2}) cannot be directly 
deduced from the C.L.T. 
\item[$\bullet$] when $d\ge3$, Theorem \ref{th2} implies that the walk is transient. 
  This result was first proved in  \cite{grimmett93a}. 
\item[$\bullet$] we comment  a little on the lower bound for the kernel 
of $X_t$ in Appendix \ref{app-lower}.  
\item[$\bullet$] D. Heicklen and C. Hoffman \cite{heicklen??a} also obtained upper estimates for 
$P_0^{\omega}[X_t=y]$ using a different method than ours but they missed the correct 
limit behaviour by a logarithmic factor. 
\end{itemize}

\subsection{Isoperimetric inequalities}

It is known that the large time behaviour of a Markov chain on a graph 
is related to the geometry at infinity of the graph. 
For instance, isoperimetric inequalities or Nash inequalities imply  
estimates on the heat kernel 
decay, {\it i.e.} upper bounds on $\sup_{x,y} P_x[X_t=y]$ (cf. Coulhon 
\cite{coulhon??a} or Pittet and Saloff-Coste \cite{pittet??a}). 

 We also use isoperimetric inequalities here, but, 
as we consider non-regular graphs such as percolation clusters, the 
classical approaches do not directly apply. 
In particular, since a percolation cluster contains, with probability one, 
an arbitrarily long linear piece, it is easy to see that 

 Q.a.s., as $t$ tends to $+\infty$,  
\begin{equation}
\label{ineq-perco} 
  \sup_{x,y\in\CC}P^\omega_x[X_t=y]
  =   
  O(t^{-1/2})\,.
\end{equation}

In the note \cite{mathieu01a}, we sketched the proof of Theorems \ref{th1} and \ref{th2} using 
local isoperimetric inequalities. 
Our proof here will be slightly different. Instead of estimating the decay of the kernel of the 
random walk $X$ killed when leaving a big box, as in the note, we shall rather estimate the 
probability transitions of the random walk restricted to a finite box. The computation becomes 
a little heavier but we obtain estimates on the mixing time for the walk confined to a box 
that have their own interest.

Let us define $\CC^{n}$ to be the connected component of $\CC\cap [-n,n]^d$ that contains the 
origin. Note that, by this definition, a point is in $\CC^n$ if and only if it belongs to 
$\CC$ and can be reached from the origin by an open path contained in $[-n,n]^d$. 
 
Let $\eps >d$ and define the isoperimetric constant 
\begin{displaymath}
I_{\eps}(\CC^n)
=
\inf_{A \subset \CC^n\,,\# A \leq \frac{\# \CC^n}{2}} 
\frac{\# (\partial_{\CC^n} A)}{(\# A)^{\frac{\eps-1}{\eps}}}\,,
\end{displaymath}
where $\partial_{\CC^n} A$ is the boundary of the set $A$ in $\CC^n$, {\it i.e.} the set of nearest neighbor
points $x \in \CC^n$ and $y \in \CC^n$ such that $\omega(x,y)=1$ and 
with either $x \in A$ and $y \not\in A$ or $x \not\in A$ and $y\in A$.
In Sections \ref{isop-2} and \ref{isop-d}, we prove that: 

  under the assumptions of Theorem \ref{th1} or \ref{th2}, 
for some constant $\beta$ that depends only on $p$ and $d$, $Q$.a.s. on the set 
$\# \CC=+\infty$, for large enough $n$, one has the inequality:
\begin{equation}
\label{ineg-isop}
I_{\eps(n)}(\CC^n) \geq \frac{\beta}{n^{1-\frac{d}{\eps(n)}}}\,,
\end{equation} 
where $\eps(n)= d+2\,d\,\frac{\log \log n}{\log n}$.


\subsection{The random walk on a finite box}

Let $(X^n_t,\ t \geq 0)$ be the random walk $X$ restricted to the set $\CC^n$. 
The definition of $X^n$ is the same as for $X$ except that jumps outside $\CC^n$ are 
now forbidden: 
the random walker waits for an exponential time of parameter $1$ and
 then chooses, uniformly at random,  
one of its neighbors in $\Z^d$, say $y$. If 
$\omega(x,y)=1$ and $y\in\CC^n$, then the walker moves to $y$; otherwise it stays still. 
Thus $X^n_t$ defines 
a Markov chain on $\CC^n$ which is reversible with respect to the counting measure on $\CC^n$. 

It follows from general considerations on finite Markov chains, see Saloff-Coste 
\cite{saloff96a}, that the isoperimetric inequality (\ref{ineg-isop}) yields different estimates 
on the kernel of $X^n$. 
Indeed (\ref{ineg-isop}) implies the following Nash inequality:

  under the assumptions of Theorem \ref{th1} or \ref{th2}, 
there exists a constant $\beta$ such that $Q$ a.s. on the set 
$\# \CC=+\infty$,  
$\exists n_0(\omega)$ and $\forall n \geq n_0(\omega)$, 
\begin{equation}
\label{nash-inequality}
 \Var(g)^{1+\frac{2}{\eps(n)}}
  \leq 
  \frac{8}{\beta^2}\,
  n^{2(1-\frac{d}{\eps(n)})}\,\EE^n(g,\,g)\,\|g\|_1^{4/\eps(n)}\,,
\end{equation}
where $\EE^n(\cdot,\,\cdot)$ is the Dirichlet form of the Markov chain $X^n$. 
The variance and the $L_1$ norms are computed
with respect to the uniform probability measure on $\CC^n$.  
Inequality (\ref{nash-inequality}) is a direct application of  
Theorem 3.3.11 of Saloff-Coste \cite{saloff96a}. Besides (see Theorem 2.3.1 of 
\cite{saloff96a}), the Nash inequality (\ref{nash-inequality}) implies  estimates 
on the transition probability:
\begin{equation}
\label{estim-reflected} 
\sup_{x,y\in\CC^n}
\left| \frac{1}{\# \CC^n} - P^\omega_x[X^n_t=y] \right|
\leq 
\left(\frac{4\,\eps(n)}{\beta^2}\right)^{\frac{\eps(n)}{2}}
\,\frac{n^{\eps(n)-d}}{t^\frac{\eps(n)}{2}}\,.
\end{equation}

Another consequence of the isoperimetric inequality is a lower bound on the spectral gap: 
let $\lambda^n$ be the lowest non zero eigenvalue of the discrete Laplacian on 
$\CC ^n$. 
From Cheeger's inequality and the isoperimetric inequality (\ref{ineg-isop}),
we deduce (see Lemma 3.3.7 of \cite{saloff96a}):
\begin{theorem}
Under the assumptions of Theorem \ref{th1} or \ref{th2},  
there exists a constant $\beta$ such that $Q$.a.s. on the set 
$\# \CC=+\infty$, for large enough $n$
\begin{displaymath}
 \lambda^n
 \geq 
 \frac{\beta}{n^2}\,.
\end{displaymath} 
\end{theorem}
See also I. Benjamini and E. Mossel \cite{benjamini00a} 
\footnote{ I. Benjamini asked us to mention that there is a gap in the 
renormalization argument given in \cite{benjamini00a}. Compare with 
parts \ref{another} and \ref{renormalization} here.}.

\section{The isoperimetric inequality when $d=2$}
\label{isop-2}

We consider the site percolation model in $\Z^2$, described in Section \ref{intro-2d}.
For $x, y \in \Z^2$, we write: $x \sim y$ if $x$ and $y$ are neighbors, and 
$\E=\{(x,y) \in \CC \times \CC\,,\,x \sim y\}$.

\subsection{A preliminary result}
\label{percolation-section}

We consider the box  
$B_{m,\,n}\eqdef([0,\,m] \times [0,\,n]) \cap \Z^2$.

\begin{figure}[h]
\centering
\includegraphics[scale=0.3]{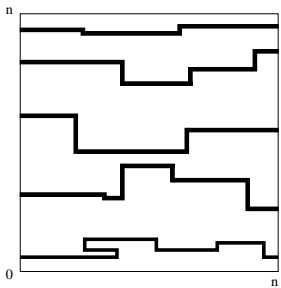} 
\caption{Example of horizontal channels}
\label{channels}
\end{figure}

\begin{definition} 
A {\it horizontal} (resp. {\it vertical}) {\it channel} of $B_{m,\,n}$  (see 
Figure \ref{channels}) 
is a path $(v_0,\,e_1,\hdots, e_n,v_n)\,,$ $v_i \in \Z^2$ and $e_i \in \E$
for all $0 \leq i \leq n\,,$ such that
\begin{itemize} 
\item 
$(v_1,\,e_1,\hdots, e_{n-1},v_{n-1})$ 
 is contained in the interior of $B_{m,\,n}\,,$ 
\item 
$v_0 \in \{0\} \times [0,\,n]$ (resp. $v_0 \in [0,\,m]\times
\{0\}$)\,, 
\item 
$v_n \in \{m\} \times [0,\,n]$ (resp. $v_n \in [0,\,m]\times \{n\}$)\,.
\end{itemize}

We say that two channels are {\em disjoint} if they do not have any 
vertex in common. 
\end{definition}  

Let $ N(m,\,n)$ be the maximal number of disjoint open  
horizontal channels in $B_{m,\,n}.$ 

\begin{theorem}[{\sf Kesten}\cite{kesten82a}, Theorem 11.1] 
\label{kesten-thm} 

Let $p > p_c$.  
\\ 
For some constant $c(p)>0$, and for some universal constants  
$0<c_1,\,c_2,\,\delta <+\infty$, one has:  
\begin{displaymath} 
 Q(N(m,\,n) \geq c(p)\,n) 
 \geq 
 1-c_1\,(m+1) \exp(-c_2\,(p-p_c)^\delta n)\,. 
\end{displaymath} 

\end{theorem} 

\subsubsection*{Construction of the {\em Kesten grid}}
\label{kesten-grid}

Theorem \ref{kesten-thm} gives us an information about the
``geometrical structure'' of $\CC^n$: 
the number of non--intersecting horizontal and vertical channels, which cross the box and 
belong to the infinite cluster, is proportional to the size of the box. They form a grid, that we call the 
{\em Kesten grid}. Let us construct and consider the following particular 
Kesten grid: we divide the box $[-n,n]^2$ in horizontal strips of width 
$C\,\log n$, with the constant $C$ large enough so that the 
expression $\ c_1\,(n+1) \exp(-c_2\,(p-p_c)^\delta\, C\,\log n)\,$ 
is summable (and therefore one may apply  Borel-Cantelli lemma).
Then it follows from Theorem \ref{kesten-thm}, applied in each strip, that there 
is a number proportional to $C\,\log n$ of horizontal channels. 
We do the same construction for vertical channels as well.

\begin{figure}[h]
\centering
\includegraphics[scale=0.3]{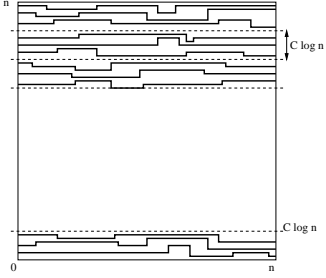} 
\caption{In each strip of width $C\,\log n$, we consider the Kesten's channels}
\label{kestenbox}
\end{figure}

The Kesten grid so constructed is regularly spread all over the
box $[-n,n]^2$ (cf. Figure \ref{kestenbox}).

\subsection{Proof of the isoperimetric inequality}

Choose $\eps>2$, which depends on $n$ and satisfies the condition:
\begin{equation}
\label{cond-epsilon} 
\liminf_{n\rightarrow+\infty} 
 (\eps -2) \,\frac{\log n}{\log \log n}
 \geq 4\,.
\end{equation} 
Condition (\ref{cond-epsilon}) is in particular satisfied if 
$\eps=2+4\frac{\log \log n}{\log n}$. 
We want to prove the inequality (\ref{ineg-isop}) in the case $d=2$, {\it i.e.}:
there exists some constant $\beta$ that depends only on $p$ such that $Q$.a.s. on the set 
$\# \CC=+\infty$, for large enough $n$, one has the inequality:
\begin{displaymath}
I_{\eps(n)}(\CC^n) 
\geq 
\frac{\beta}{n^{\frac{\eps-2}{\eps}}}\,. 
\end{displaymath}

We first remark that, without loss of generality, in the definition of the isoperimetric
constant $I_{\eps}(\CC^n)$, we can take the infimum only  on connected subsets of 
$\CC^n$.

Let $A$ be a finite connected subset of $\CC^n$. 
We denote by $N$ the cardinal of  $A$ ($N \eqdef \# A$). 
We look for a lower bound on the cardinal of $\partial_{\CC^n} A$.
The following result will be very useful:

$\exists \alpha(p) >0$ such that, $Q$-a.s on the event  
$\# \CC=+\infty$,   
\begin{equation}
\label{cardinal-cluster}
 \alpha(p) \leq 
  \Liminf_{n \to \infty} \frac{\# \CC^n}{n^2}
 \leq 1\,.
\end{equation} 
Indeed, the first inequality (left) is a consequence of Theorem
\ref{kesten-thm}, and the second one comes from $\# \CC^n \leq n^2$.

\medskip

We introduce several classes of sets $A$ and derive a lower bound on 
$\frac{\#(\partial_{\CC^n} A)}{(\# A)^{\frac{\eps-1}{\eps}}}$ in each of them.
 
\begin{enumerate}
\item Let 
$\AAA_0 = \{ A \subset \CC^n, A \mbox{ connected }, \ N \leq n^{\frac{\eps -2}{\eps -1}}\}$.
 
Obviously $\#\partial A \geq 1$, so we have:
\begin{equation}
\label{cas_0}
\min_{A \in \AAA_0}
  \frac{\#(\partial_{\CC^n} A)}{(\# A)^{\frac{\eps-1}{\eps}}}
 \geq
 \frac{1}{N^{\frac{\eps-1}{\eps}}}
 \geq
 \frac{1}{n^{\frac{\eps-2}{\eps}}}\,.
\end{equation}

\item Consider now sets $A$ such that  
$n^{\frac{\eps -2}{\eps -1}} < N \leq \frac{\# \CC^n}{2}$ .
\end{enumerate}

We consider the Kesten grid, constructed in Section \ref{percolation-section}.
Two cases appear naturally: 
{\it (1)} $A$ contains either no horizontal channel or no vertical one, or
{\it (2)} $A$ contains at least one horizontal and one vertical channel.

If $R_A$ denotes the smallest rectangle which contains the set $A$, then  
$R_A$ is strictly included in $[-n,n]^2$ in the first case, and 
$R_A=[-n,n]^2$ in the second one.

We call $\diam(A)$ the length of the longest side of $R_A$.

\subsubsection*{{\it (1)} $A$ contains either no horizontal channel
                or no vertical channel}

\begin{itemize}

    \item[$\bullet$] Let 
$\AAA_1 = 
\{ A \subset \CC^n, A \mbox{ connected }, 
   n^{\frac{\eps -2}{\eps -1}} < N \leq \frac{\# \CC^n}{2},\ \diam(A) <n\}$.

We know that every channel which 
intersects both the set $A$ and its complementary $\CC^n \setminus A$ contains
at least one element of $\partial_{\CC^n} A$.

Let us consider the $N_A$ channels which intersect $R_A$ (we consider 
the horizontal channels if $\diam(A)$ is the vertical side, the vertical 
channels otherwise). Each of these channels brings at least a contribution 
of $1$ in the cardinal of $A$, so: $\ \# \partial_{\CC^n} A \geq N_A$.

\begin{figure}[h]
\centering
\includegraphics[scale=0.3]{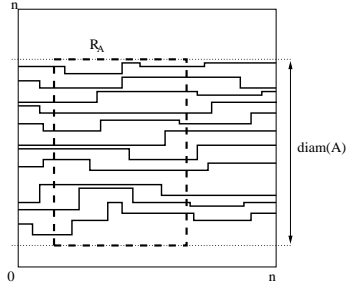}
\caption{The number of horizontal channels which intersect $R_A$ is 
         proportional to $\diam(A)$}
\label{cas2}
\end{figure}

Thanks to condition (\ref{cond-epsilon}),  
\begin{displaymath}
 \diam(A) \geq \sqrt{\# A} 
  = \sqrt{N} \geq n^{\frac{\eps -2}{2(\eps -1)}}\,,
\end{displaymath} 
and this is larger than  $C \log n$, for any constant $C$. 
Therefore  Kesten's theorem can be applied, {\it i.e.} \\
$\ N_A \geq c_1(p)\,\diam(A)\ $  and
$\quad\# \partial_{\CC^n} A \geq c_1(p)\,\diam(A) \geq c_1(p)\,\sqrt{N}$. 

Thus, for any connected set 
$A \subset \CC^n,\, n^{\frac{\eps -2}{\eps -1}} < N \leq
\frac{\# \CC^n}{2}$, 
\begin{displaymath}
  \frac{\#(\partial_{\CC^n} A)}{(\# A)^{\frac{\eps-1}{\eps}}}
 =
 \frac{\# (\partial_{\CC^n} A)}{N^{\frac{\eps-1}{\eps}}}
 \geq
 \frac{c_1(p)\,\sqrt{N}}{N^{\frac{\eps-1}{\eps}}}\,.
\end{displaymath}
This last quantity  reaches its minimum for $N=\frac{\# \CC^n}{2}$ and:
\begin{equation}
\label{cas_1}
\min_{A \in \AAA_1} 
  \frac{\#(\partial_{\CC^n} A)}{(\# A)^{\frac{\eps-1}{\eps}}}
 \geq
2^{\frac{\eps-2}{2\eps}}\ \frac{c_1(p)}{n^{\frac{\eps-2}{\eps}}}\,.
\end{equation}

   \item[$\bullet$] Let 
 $\AAA_2 = 
 \{ A \subset \CC^n, A \mbox{ connected }, 
    n^{\frac{\eps -2}{\eps -1}} < N \leq \frac{\# \CC^n}{2},\ \diam(A) 
    =n\}.$

 We use the same reasoning as in the previous case, but here, all the 
horizontal (or vertical) channels of the Kesten grid intersect $R_A$, 
since $\diam(A)=n$ (and so, with Theorem \ref{kesten-thm}, 
$N_A \geq c(p)\,n$).
Hence, $\quad \# \partial_{\CC^n} A \geq N_A \geq c(p)\,n$, and
\begin{displaymath}
  \frac{\#(\partial_{\CC^n} A)}{(\# A)^{\frac{\eps-1}{\eps}}}
 =
 \frac{\# \partial_{\CC^n} A}{N^{\frac{\eps-1}{\eps}}}
 \geq
 \frac{c(p)\,n}{N^{\frac{\eps-1}{\eps}}}\,.
\end{displaymath}
The minimum is reached for  $N=\frac{\# \CC^n}{2}$:
\begin{equation}
\label{cas_2}
\min_{A \in \AAA_2}
      \frac{\#(\partial_{\CC^n} A)}{(\# A)^{\frac{\eps-1}{\eps}}}
  \geq
  2^{\frac{\eps-1}{\eps}}\ \frac{c(p)}{n^{\frac{\eps-2}{\eps}}}\,.
\end{equation}

\end{itemize}

\subsubsection*{{\it (2)} $A$ contains at least one horizontal
                {\sc and} one vertical channel}

It implies that $R_A = [-n,n]^2$.

In this case, one has to take into account that 
Kesten's channels which are completely included in the set
$A$ do not contribute in $\# \partial_{\CC^n} A$. 
Either there is, at least in one direction, a non--negligeable 
proportion of Kesten's channels which are not completely included in $A$, or  
in both directions, almost all the channels are 
completely included in $A$. We consider separately these two possibilities. 

\bigskip

The Kesten grid is constituted of 
$N(n)\geq c(p)\,n$ horizontal and vertical channels. 
We distinguish channels which are completely included in $A$ (i.e 
which have an empty intersection with $\CC^n \setminus A$).
Let:
\begin{displaymath}
 \delta_H 
 \eqdef 
 \frac{1}{n}\,\# \{ \mbox{ horizontal channels } \subset A \}\,,
 \qquad 
 \delta_V 
 \eqdef 
 \frac{1}{n}\,\# \{ \mbox{ vertical channels }  \subset A \}\,.
\end{displaymath}
For some $\delta$ to be chosen later in the interval $]0,c(p)[$, we define:

\begin{itemize}
 \item  $\AAA_3 = 
 \{ A \subset \CC^n, A \mbox{ connected }, 
    n^{\frac{\eps -2}{\eps -1}} < N \leq \frac{\# \CC^n}{2},\ R_A = 
    [-n,n]^2,\,\\
    \delta_H < \delta \mbox{ or } \delta_V < \delta\}.$

 At least in one direction, there are less than $\delta\,n$ channels
which are completely included in $A$. As the total number of channels
in one direction is $c(p)\,n$, we have more than $(c(p)-\delta)\,n$
channels which are {\sc not} completely included in $A$, and therefore 
intersect $\CC^n \setminus A$. Thus: 
$\ \# \partial_{\CC^n} A \geq (c(p)-\delta)\,n\,,$ and
\begin{displaymath}
  \frac{\#(\partial_{\CC^n} A)}{(\# A)^{\frac{\eps-1}{\eps}}}
 \geq
 \frac{(c(p)-\delta)\,n}{N^{\frac{\eps-1}{\eps}}}\,.
\end{displaymath}
The minimum is reached for  $N=\frac{\# \CC^n}{2}$:
\begin{equation}
\label{cas_3}
\min_{A \in \AAA_3}
      \frac{\#(\partial_{\CC^n} A)}{(\# A)^{\frac{\eps-1}{\eps}}}
  \geq
  2^{\frac{\eps-1}{\eps}}\ \frac{(c(p)-\delta)}{n^{\frac{\eps-2}{\eps}}}\,.
\end{equation}

 \item Let 
\begin{eqnarray*}
 \AAA_4 &=& 
 \{ A \subset \CC^n, A \mbox{ connected }, 
    n^{\frac{\eps -2}{\eps -1}} < N \leq \frac{\# \CC^n}{2},\\
    && \qquad\qquad\qquad\qquad R_A = [-n,n]^2,\, 
    \delta_H > \delta \mbox{ and } \delta_V > \delta\}.
\end{eqnarray*} 
\end{itemize}

\begin{figure}[h]
\centering
\includegraphics[scale=0.3]{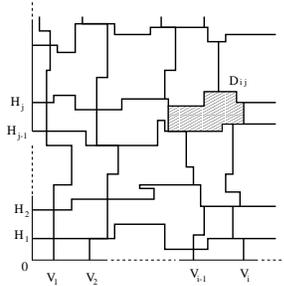}
\caption{The domain $D_{i,\,j}$ is included in the grey part}
\label{domain}
\end{figure}

Let us consider the sub--grid of the Kesten grid, constituted only with the 
channels included in $A$.
These channels divide $[-n,n]^2$ in different non-overlapping regions,
which will be called {\em domains}.

If we enumerate the horizontal channels of the sub-grid by 
$H_1,\,H_2,\,\hdots$, and the vertical ones by $V_1,\,V_2,\,\hdots$, the 
domain $D_{i,\,j}$ is constituted of all the points of $\CC^n$ included 
(strictly) between $V_{i-1}$, $V_i$, 
$H_{j-1}$ and $H_j$ (with: $H_0=[0,\,n]\times{0}$ and 
 $V_0= {0} \times [0,\,n]$). See Figure \ref{domain}.

Each of these domains which intersects $\CC^n \setminus A$ increases by one
unit the cardinal of $\# \partial_{\CC^n} A$.
Thus, one way to find a lower bound for $\# \partial_{\CC^n} A$ is to count the
number of such domains.

\begin{figure}[h]
\centering
\includegraphics[scale=0.3]{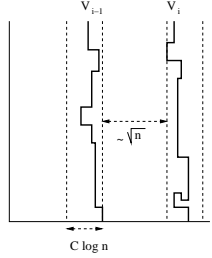}
\caption{Fluctuations of channels are not greater than $C\,\log n$}
\label{dist-channels}
\end{figure}

We classify  domains into two categories:  {\em big} 
and  {\em small} ones. To do so, we need to define the 
{\em distance between two channels}. 
(Channels are irregular, but thanks to the construction of 
the channels in strips (cf Figure \ref{kestenbox}), their
fluctuations are not greater than $C\,\log n$. See Figure
\ref{dist-channels}.) 
\begin{definition}
Let $C_1$ and $C_2$ be two channels in the box $[-n,n]^2$ (both in
the same direction). 
We define the distance between $C_1$ and $C_2$ 
by:
\begin{displaymath}
 d(C_1,\,C_2)
 \eqdef
 \inf \{d(x,\,y)\,;\ x \in C_1,\ y \in C_2 \}
\end{displaymath} 
where $d:\,\Z^2\times \Z^2 \rightarrow \R_+$ is the usual Euclidean
distance.
\end{definition}

We will consider that the domain $D_{i,\,j}$ is {\it big} if the distance
$d(V_{i-1},\,V_i)$ or $d(H_{j-1},\,H_j)$ is greater than or equal to  
$\sqrt{n}$. The other domains are {\it small}. 
Let $B$ be the set of  {\em big domains}. 

Let us evaluate the total volume of big domains.
First, we consider the set $B_V$ constitued of all the domains
$D_{i,\,j}$ such that $d(V_{i-1},\,V_i)\geq \sqrt{n}$.
Clearly, the volume of $B_V$ is smaller than the volume between all the
channels $V_{i-1}$ and $V_i$ such that $d(V_{i-1},\,V_i)\geq
\sqrt{n}$. 
The domain delimited by the channels $V_{i-1}$ and $V_i$ 
contains  a rectangle of size $\ n \times d(V_{i-1},\,V_i)$, and 
it is strictly included in the rectangle 
$\,n \times [d(V_{i-1},\,V_i)+2\,C\,\log n]$, because the fluctuations
of the channels are smaller than $C\,\log n$ (see Figure
\ref{dist-channels}). Therefore we have the following bound:
\begin{eqnarray*}
 \# B_V
 &\leq &
 n \left[
 \sum_{\{i;\,d(V_{i-1},\,V_i)\geq \sqrt{n}\}}
 d(V_{i-1},\,V_i) + 2\,C \log n
 \right]
\\
 &\leq&
 n\,\sum_{\{i;\,d(V_{i-1},\,V_i)\geq \sqrt{n}\}}d(V_{i-1},\,V_i)
 +
 2\,n\,\sqrt{n}\,C \log n \,.
\end{eqnarray*} 
On the scale $\sqrt{n}$, we may apply  Theorem
\ref{kesten-thm}: if $d(V_{i-1},\,V_i)\geq \sqrt{n}$, 
there is $c(p)\,d(V_{i-1},\,V_i)$ Kesten channels between the channels
$V_{i-1}$ and $V_i$.
These channels belong to the $(c(p)-\delta_H)\,n$ ones  which have not been 
considered in the sub--grid ({\it i.e.} which are not completely included in $A$).
Therefore,
\begin{displaymath}
 \sum_{\{i;\,d(V_{i-1},\,V_i)\geq \sqrt{n}\}}
 c(p)\,d(V_{i-1},\,V_i)
 \leq (c(p)-\delta_H)\,n
 \leq (c(p)-\delta)\,n
\end{displaymath}
and 
\begin{displaymath}
\sum_{\{i;\,d(V_{i-1},\,V_i)\geq \sqrt{n}\}}
 d(V_{i-1},\,V_i)
 \leq 
\frac{(c(p)-\delta)}{c(p)}\,n\,.
\end{displaymath} 
Therefore
\begin{displaymath}
 \# B_V
 \leq 
 \frac{(c(p)-\delta)}{c(p)}\,n^2
 +2\,n^{3/2}\,C\,\log n\,.
\end{displaymath}
We find the same upper bound for the volume of the
set $B_H$, constitued of all the domains
$D_{i,\,j}$ such that $d(H_{j-1},\,H_j)\geq \sqrt{n}$.
So,
\begin{displaymath}
 \# B
 \leq
 2\,\frac{(c(p)-\delta)}{c(p)}\,n^2
 +4\,n^{3/2}\,C\,\log n\,.
\end{displaymath} 
It is now sufficient to take $\delta$ such that 
$\ \displaystyle\frac{(c(p)-\delta)}{c(p)} 
  \leq \frac{\# \CC^n}{9\,n^2}$
to see that, for large enough $n$,
\begin{displaymath}
 \# B \leq \frac 3 9 \, \# \CC^n\,.
\end{displaymath}
Therefore 
\begin{displaymath}
 \#(B\cap\CC^n) \leq \frac 3 9 \, \# \CC^n\,.
\end{displaymath}
Call $S$ the union of small domains. 
Thus $S\cap (\CC^n\setminus A)$  covers an area of at least $\frac 6 9 \# \CC^n$. 
We have assumed  that $\# (\CC^n \setminus A) \geq \frac{\# \CC^n}{2}$. 
Therefore 
\begin{displaymath}
 \#(S\cap (\CC^n\setminus A)) \geq \frac 1 6 \, \# \CC^n\,.
\end{displaymath}
Besides, by definition, each small domain
has  volume smaller than $n$.
Thus,  
\begin{displaymath}
\# \left\{ D_{i,\,j}\,:\, 
  D_{i,\,j} \mbox{ small, and } D_{i,\,j} \cap (\CC^n \setminus A) \neq \emptyset
    \right\} 
 \geq \,\frac{\# \CC^n}{6n}\,.
\end{displaymath}
Since any such domain $D_{i,j}$ contributes by at least one in $\#\partial_{\CC^n} A$,  
it follows that: $\ \# \partial_{\CC^n} A \geq \frac{\# \CC^n}{6n}\,,$ and
\begin{displaymath}
  \frac{\#(\partial_{\CC^n} A)}{(\# A)^{\frac{\eps-1}{\eps}}}
 \geq
 \frac{\# \CC^n}{6n\,N^{\frac{\eps-1}{\eps}}}\,.
\end{displaymath}
We conclude that 
\begin{equation}
\label{cas_4}
\min_{A \in \AAA_4} 
      \frac{\#(\partial_{\CC^n} A)}{(\# A)^{\frac{\eps-1}{\eps}}}
  \geq
  \frac {2^{\frac{\eps-1}{\eps}}} {6n^{\frac{\eps-2}{\eps}}}\,.
\end{equation}

To conclude, we just gather the inequalities 
(\ref{cas_0}), (\ref{cas_1}), (\ref{cas_2}), (\ref{cas_3}) and 
(\ref{cas_4}).

\carre

\section{The isoperimetric inequality when $d\geq 2$}
\label{isop-d}

We consider the bond percolation model in $\Z^d$, described in section \ref{intro-d}.
For $x, y \in \Z^d$, we write: $x \sim y$ if $x$ and $y$ are neighbors in $\Z^d$, and 
$\E_d=\{(x,y) \in \Z^d \times \Z^d\,,\,x \sim y\}$.
The application $\omega:\E_d \rightarrow [0,1]$ is called a configuration, and the random 
variables $(\omega(e),\,e \in \E_d)$ are Bernouilli$(p)$ independent variables. 

Let $\BB^n=[-n,n]^d$. 
We have already defined the isoperimetric constant:
\begin{equation}
\label{cst-isop2}
I_{\eps}(\CC^n)
=
\inf_{A \subset \CC^n\,,\# A \leq \frac{\# \CC^n}{2}} 
\frac{\# (\partial_{\CC^n} A)}{(\# A)^{\frac{\eps-1}{\eps}}}\,,
\end{equation}
with $\partial_{\CC^n} A = \left\{ (x,y) \in \E_d\,, \omega(x,y)=1\,:\,\left|
      \begin{array}{l} x \in A \cap \CC^n \\ y \in A^C\cap \CC^n \end{array}\right.
      \mbox{ or } \left|
      \begin{array}{l} x \in A^C \cap \CC^n \\ y \in A\cap \CC^n \end{array}\right.
      \right\}$ 
and $A^C$ the complement of $A$ in $\Z^d$. 

It will be convenient to also introduce the isoperimetric constants:
\begin{equation}
\label{cst-isop3}
I_{\eps}^{(\alpha)}(\CC^n)
=
\inf_{A \subset \CC^n\,,\# A \leq (1-\alpha) \# \CC^n }
\frac{\# (\partial_{\CC^n} A)}{(\# A)^{\frac{\eps-1}{\eps}}}\,,
\end{equation} 
where $\alpha\in]0,1/2]$. Note that $I_\eps(\CC^n)=I_\eps^{(\frac 12)}(\CC^n)$. 

We claim that: 

for $\alpha\in]0,\frac 12]$ and $p>p_c$, there exists a constant $\beta$ such that, 
Q.a.s. on the set $\#\CC=\infty$, and for large enough $n$, one has 
\begin{equation}
\label{*} 
I_{\eps(n)}^{(\alpha)}(\CC^n)\ge \beta n^{\frac d{\eps(n)}-1}\,,
\end{equation}
where $\eps(n)=d+2d\frac{\log\log n}{\log n}$. 

The proof goes in three steps.

\subsection{Geometric arguments}
\label{geometric}

We first prove that (\ref{*}) is equivalent to (\ref{*2}). This proof 
is based on general arguments that would actually work on any graph. Next, we prove that 
(\ref{*2}) is implied by (\ref{*3}). This second step relies on specific properties of the 
underlying graph $\Z^d$. Although both proofs are rather classical, we prefered to give 
some details for the reader's convenience. 

Choose $\alpha\in]0,\frac 12]$. 
Define 
\begin{displaymath}
\bar{I}_{\eps}(\CC^n)
=
\inf_{A \subset \CC^n} 
\frac{Q^n(\partial_{\CC^n} A)}{(\pi^n(A)\, \pi^n(\CC^n\setminus A))^{\frac{\eps-1}{\eps}}}\,,
\end{displaymath}
where $\pi^n$ is the uniform probability on $\CC^n$ and, \\
$\forall B \subset \E_d\,,$
$Q^n(B)=\frac{1}{2\,d\,\# \CC^n} \# \{ e=(x,y) \in B, x,y \in \CC^n, \omega(x,y)=1 \} $. 
Then
\begin{displaymath}
\bar{I}_{\eps}(\CC^n) \leq \alpha^{\frac 1 \eps -1} 
\inf_{A \subset \CC^n\,,\# A \leq (1-\alpha)\# \CC^n} 
\frac{Q^n(\partial_{\CC^n} A)}{(\pi^n(A))^{\frac{\eps-1}{\eps}}}
\leq
\alpha^{\frac 1 \eps -1}  \,{\bar I}_{\eps}(\CC^n)\,.
\end{displaymath} 
The first inequality follows from the fact that, if $\ \# A \leq (1-\alpha)\# \CC^n$, 
then \\ $\pi^n(\CC^n\setminus A)\geq \alpha$. To get the second inequality note that, 
since $\alpha\in]0,\frac 12]$, for any 
set $A\subset \CC^n$, either $\# A \leq (1-\alpha)\# \CC^n$ or 
$\#(\CC^n\setminus A)  \leq (1-\alpha)\# \CC^n$.  
Also, since $\pi^n(A)=\frac{\# A}{\# \CC^n}$, 
\begin{displaymath}
\inf_{A \subset \CC^n\,,\# A \leq (1-\alpha)\# \CC^n} 
\frac{Q^n(\partial_{\CC^n} A)}{(\pi^n(A))^{\frac{\eps-1}{\eps}}}
= \frac { (\#\CC^n)^{-\frac 1\eps}}{2d}  I_\eps^{(\alpha)} (\CC^n)\,.
\end{displaymath}
In particular, it is equivalent to prove (\ref{*}) for all values of $\alpha\in]0,\frac 12]$ or 
for some given value of $\alpha\in]0,\frac 12]$. 

Moreover, in the computation of $\bar{I}_{\eps}(\CC^n)$, the infimum is reached for sets $A$
such that $A$ and  $\CC^n \setminus A$ are connected in $\CC^n$. 
Indeed,  write $A= \cup_{i} A_i$, where each set $A_i$ is connected in $\CC^n$.\\
 
Let: $ \gamma = \displaystyle \inf_i 
  \frac{Q^n (\partial_{\CC^n} A_i)}
       {(\pi^n(A_i) \,\pi^n(\CC^n \setminus A_i))^{\frac{\eps-1}{\eps}}}$. Then,
\begin{eqnarray*}
 Q^n(\partial_{\CC^n} A)
 =
 \sum_i Q^n(\partial_{\CC^n} A_i)
 &\geq&
 \gamma\, \sum_i (\pi^n(A_i) \,\pi^n(\CC^n \setminus A_i))^{\frac{\eps-1}{\eps}}
 \\
 &=&
 \gamma\, \sum_i (\pi^n(A_i) \,(1-\pi^n(A_i))^{\frac{\eps-1}{\eps}}
 \\
 &\geq&
 \gamma\, (\sum_i \pi^n(A_i) \,(1-\pi^n(A_i)))^{\frac{\eps-1}{\eps}}
 \\
 &=&
 \gamma \,(\pi^n(A)-\sum_i \pi^n(A_i)^2)^{\frac{\eps-1}{\eps}} \,.
\end{eqnarray*}
From  $ \sum_i \pi^n(A_i)^2 \leq \left( \sum_i \pi^n(A_i)\right)^2=(\pi^n(A))^2$, we get:\\
$Q^n(\partial_{\CC^n} A) \geq \gamma\,(\pi^n(A)(1-\pi^n(A)))^{\frac{\eps-1}{\eps}}\,$, 
{\it i.e.} 
\begin{displaymath}
\frac{Q^n (\partial_{\CC^n} A)}
       {(\pi^n(A) \,\pi^n(\CC^n \setminus A))^{\frac{\eps-1}{\eps}}}
       \geq
 \gamma=
\inf_i \frac{Q^n (\partial_{\CC^n} A_i)}
       {(\pi^n(A_i) \,\pi^n(\CC^n \setminus A_i))^{\frac{\eps-1}{\eps}}}\,.
\end{displaymath} 
We can apply the same argument to $\CC^n \setminus A_i$ instead of $A$ and prove that 
\begin{displaymath}
\frac{Q^n (\partial_{\CC^n} A)}
       {(\pi^n(A) \,\pi^n(\CC^n \setminus A))^{\frac{\eps-1}{\eps}}}
       \geq 
\inf_{i,j} \frac{Q^n (\partial_{\CC^n} A_{i,j})}
       {(\pi^n(A_{i,j}) \,\pi^n(\CC^n \setminus A_{i,j}))^{\frac{\eps-1}{\eps}}}\,,
\end{displaymath}
where, for all $i$, the sets $A_{i,j}$ are the connected components of $\CC^n \setminus A_i$. 
Since the sets $A_{i,j}$ are connected and such that $\CC^n \setminus A_{i,j}$  is connected, 
we have indeed proved that the infimum in the definition  of $\bar{I}_{\eps}(\CC^n)$ is reached 
for sets $A$ such that $A$ and  $\CC^n \setminus A$ are connected in $\CC^n$.    

The same kind of remark applies to $I_\eps(\CC^n)$:  \\
choose $A\subset\CC^n$ with $\#A\le(1-\alpha)\#\CC^n$. There exists a connected subset of 
$\CC^n$, say $A_0$, such that $\CC^n\setminus A_0$ is connected and 
\begin{displaymath}
\frac{Q^n (\partial_{\CC^n} A)}
       {(\pi^n(A) \,\pi^n(\CC^n \setminus A))^{\frac{\eps-1}{\eps}}}
       \ge
\frac{Q^n (\partial_{\CC^n} A_0)}
       {(\pi^n(A_0) \,\pi^n(\CC^n \setminus A_0))^{\frac{\eps-1}{\eps}}}\,.
\end{displaymath}
Among the two sets $A_0$ and $\CC^n\setminus A_0$, one has $\pi^n$ measure smaller than 
$1-\alpha$, say $\pi^n(A_0)\le 1-\alpha$. 

Since $\#A\le(1-\alpha)\#\CC^n$, we have $\pi^n(\CC^n\setminus A)\ge \alpha$ and 
\begin{displaymath}
\frac{Q^n (\partial_{\CC^n} A)}
       {(\pi^n(A))^{\frac{\eps-1}{\eps}}}
       \geq \alpha^{1-\frac 1\eps}
 \frac{Q^n (\partial_{\CC^n} A)}
       {(\pi^n(A) \,\pi^n(\CC^n \setminus A))^{\frac{\eps-1}{\eps}}}\,.       
\end{displaymath}
Besides, 
\begin{displaymath}
\frac{Q^n (\partial_{\CC^n} A_0)}
       {(\pi^n(A_0) \,\pi^n(\CC^n \setminus A_0))^{\frac{\eps-1}{\eps}}}
  \geq
\frac{Q^n (\partial_{\CC^n} A_0)}
       {(\pi^n(A_0) )^{\frac{\eps-1}{\eps}}}\,. 
\end{displaymath}
We conclude that 
\begin{displaymath}
\frac{Q^n (\partial_{\CC^n} A)}
       {(\pi^n(A)))^{\frac{\eps-1}{\eps}}}
       \geq \alpha^{1-\frac 1\eps} 
\frac{Q^n (\partial_{\CC^n} A_0)}
       {(\pi^n(A_0) )^{\frac{\eps-1}{\eps}}}\,, 
\end{displaymath}
or equivalently  
\begin{displaymath}
\frac{\# (\partial_{\CC^n} A)}{(\# A)^{\frac{\eps-1}{\eps}}}
\ge \alpha^{1-\frac 1\eps} 
\frac{\# (\partial_{\CC^n} A_0)}{(\# A_0)^{\frac{\eps-1}{\eps}}}\,. 
\end{displaymath} 
We conclude that, in the definitions of the isoperimetric constants 
$I_\eps(\CC^n)$ and $I_\eps^{(\alpha)}(\CC^n)$, we may restrict our attention to connected sets 
$A$ such that $\CC^n\setminus A$ is connected. It will affect the value of these constants by 
at most the multiplicative constant $\alpha^{1-\frac 1\eps}$ and therefore only change the value 
of constant $\beta$ in (\ref{*}).  

\medskip 

Since, for any set $A\subset\CC^n$ with $\#A\le(1-\alpha)\#\CC^n$, 
$\partial_{\CC^n} A\not=\emptyset$, we have 
\begin{displaymath}
\frac{\# (\partial_{\CC^n} A)}{(\# A)^{\frac{\eps-1}{\eps}}}
\geq
\frac{1}{(\# A)^{\frac{\eps-1}{\eps}}}
\geq 
\beta n^{\frac d{\eps} -1}\,,
\end{displaymath}
as soon as 
$\# A \leq \frac{n^{\frac{\eps-d}{\eps-1}}}{\beta^{\frac{\eps}{\eps-1}}}$.

Thus, gathering these remarks, we see that  
an equivalent formulation of (\ref{*}) is: 

for $c>0$ and $p>p_c$, there exists a constant $\beta$ such that, 
Q.a.s. on the set $\#\CC=\infty$, and for large enough $n$, one has 
\begin{eqnarray} 
\label{*2}
&&
\inf \left\{ 
\frac{\# (\partial_{\CC^n} A)}{(\# A)^{\frac{\eps(n)-1}{\eps(n)}}}\,, 
A \subset \CC^n\,,\# A \leq \frac{\# \CC^n}{2}\,, A \mbox{ connected}\,,
\right.\\
\nonumber
&&\hspace{3cm} \left. \CC^n\setminus A \mbox{ connected}\,, 
\#A\geq c n^{\frac{\eps(n)-d}{\eps(n)-1}} \right\} 
\ge \beta n^{\frac d{\eps(n)}-1}\,. 
\end{eqnarray} 

\medskip

From now on, let $A \subset \CC^n$ be such that $A$ and $\CC^n \setminus A$ are 
connected sets and $\#A\le\frac 12 \#\CC^n$.

\begin{figure}[h]
\centering
\includegraphics[scale=0.4]{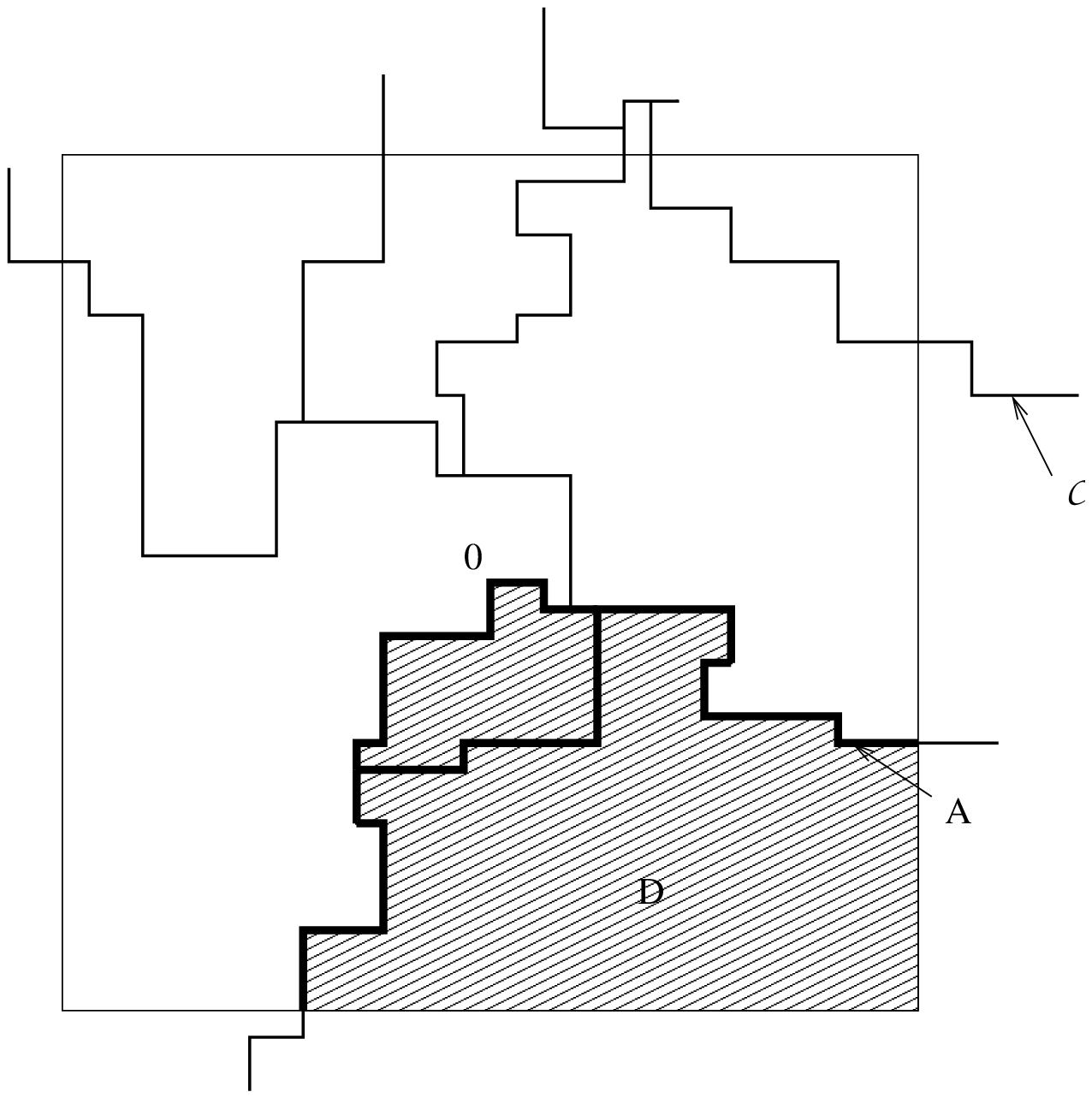} 
\caption{}
\label{connectivity}
\end{figure}

\noindent
Let $B$ be the unique $\BB^n$-connected component of $\BB^n \setminus A$ which contains 
$\CC^n \setminus A$. (Remenber that  $\BB^n=[-n,n]^d$.)  
Consider now $D$, the complementary set of $B$ in $\BB^n$~: $D=\BB^n \setminus B$ 
(cf. Figure \ref{connectivity}).
Then, we have:
\begin{enumerate}
\item 
 $D \cap \CC^n=A$, $\quad B\cap \CC^n=(\CC^n \setminus A)\,,$

\item \label{Dcon}
 $D$ is connected in $\Z^d$ (because $A$ is connected in $\Z^d$),

\item \label{DCcon}
 its complement $B=\BB^n \setminus D$ is also connected in $\Z^d$,

\item \label{front-con}
 The boundary of $A$ in $\CC^n$ satisfies: 
 
 $\partial_{\CC^n} A=\{ (x,y) \in \E_d\,|\,(x,y) \in \partial_{\BB^n} D, \
 \omega(x,y)=1\}\,.$\\
 {\bf Proof:}  Consider $(x,y) \in \partial_{\CC^n} A$, $x \in A$ and 
 $y \in (\CC^n \setminus A)$.
 Then, $x \in D$ and $y \not\in D$, so $(x,y) \in \partial_{\BB^n} D$ and $\omega(x,y)=1$.
 
  Consider $(x,y) \in \partial_{\BB^n} D$ with $x \in D$ and $\omega(x,y)=1$. Then, $x \in A$, $y
  \not\in A$ but $y \in \CC^n$ ($\CC^n$ is connected). So, $(x,y) \in \partial_{\CC^n} A$.
  \carre
 
\end{enumerate}
The two properties ({\it \ref{Dcon}}) and ({\it \ref{DCcon}}) imply that the boundary 
of $D$ in $\BB^n$:
   $\partial_{\BB^n} D=\{(x,y) \in \E_d \cap \BB^n\,|\, x \in D, y \not \in D\}$ is $*$-connected in $\Z^d$
(cf. Appendix \ref{app-con}).
Property ({\it \ref{front-con}}) implies that:  $\displaystyle 
\# \partial_{\CC^n} A = \sum_{\{e \in \partial_{\BB^n} D\}} {\bf 1}_{\omega(e)=1}\,.$

Therefore, we get
\begin{eqnarray*}
\frac {\#(\partial_{\CC^n} A)}{(\#A)^{\frac{\eps-1}\eps}} \geq
\frac{\sum_{\{e \in \partial_{\BB^n} D\}} {\bf 1}_{\omega(e)=1}}
     {(\# D)^{\frac{\eps-1}{\eps}}}\,.
\end{eqnarray*}

 Since $\# A \leq \frac{\# \CC^n}{2}$ and by property (i), 
\begin{eqnarray*}
 \# D 
 &\leq&
 \# (\BB^n \setminus \CC^n) + \# A
 \\
 &\leq& 
 \# \BB^n-\# \CC^n + \frac{1}{2} \# \CC^n
 \leq
 (2\,n+1)^d-\frac{1}{2} \# \CC^n \,.
 \end{eqnarray*}

We know (cf. Appendix \ref{app-card}) that: 

there exists a constant,
$\alpha >0$ such that $Q.a.s.$ on the set $\#\CC=+\infty$, for large $n$, 
$\# \CC^n \geq 2\alpha\,(2n+1)^d$. Therefore
\begin{displaymath}
\frac{\# D}{(2\,n+1)^d} 
 \leq
 1-\frac{1}{2} \frac{\# \CC^n }{(2\,n+1)^d} 
 \leq
 1-\alpha\,.
\end{displaymath}
Define 
\begin{displaymath}
\II^{(\alpha)}=\inf_{n} 
\inf_{D\subset\BB^n\,,\#D\le(1-\alpha)(2n+1)^d}
\frac{\#(\partial_{\BB^n}D)}{(\#D)^{\frac{d-1}d}}\,.
\end{displaymath}
The (classical) isoperimetric inequality states that $\II^{(\alpha)}>0$. 
(As a matter of fact, it is well known that $\II^{(\frac 12)}>0$. 
The general case $\alpha\in]0,\frac 12]$ follows just the same way 
we showed that (\ref{*}) holds for any $\alpha\in]0,\frac 12]$ if and only if it holds 
for some $\alpha\in]0,\frac 12]$.) 

Since, for large enough $n$, we have $\#D\leq (1-\alpha) (2n+1)^d$, we also have 
\begin{displaymath}
(\#D)^{\frac{\eps-1}\eps}
\leq \frac {\#(\partial_{\BB^n}D)}{\II^{(\alpha)}} (\# D)^{\frac{\eps-d}{\eps d}} 
\leq \frac {\#(\partial_{\BB^n}D)}{\II^{(\alpha)}} (2n+1)^{\frac{\eps-d}\eps}\,.
\end{displaymath} 

Thus, we have obtained the inequality: 
\begin{displaymath}
\frac {\#(\partial_{\CC^n} A)}{(\#A)^{\frac{\eps-1}\eps}} \geq
\frac{\sum_{\{e \in \partial_{\BB^n} D\}} {\bf 1}_{\omega(e)=1}} {\#(\partial_{\BB^n}D)} 
     \frac{\II^{(\alpha)}}{(2n+1)^{\frac{\eps-d}\eps}}\,.
\end{displaymath}

Finally note that, if $\ \#A\geq c n^{\frac {\eps-d}{\eps-1}}$, then 
$\ \#D\geq c n^{\frac {\eps-d}{\eps-1}}$ and therefore\\
$\#(\partial_{\BB^n}D)\geq c n^{\frac {(\eps-d)(d-1)}{(\eps-1)d}}$ for some other value of $c$. 

We use the notation $F=\partial_{\BB^n}D$. 
Let $\E^n$ be the set of edges in $\E_d$ with both end points in $\BB^n$. 
Since $D$ and $\BB^n\setminus D$ are connected in $\BB^n$, then $F$ is $*$-connected, 
and we conclude that, in order to prove (\ref{*}), it is sufficient to check that: 

for $c>0$ and $p>p_c$, there exists a constant $\beta$ such that, 
Q.a.s. on the set $\#\CC=\infty$, and for large enough $n$, one has 
\begin{displaymath}
\inf \{ 
\frac{\sum_{\{e \in  F\}} {\bf 1}_{\omega(e)=1}}
     {\#F}\,, 
F \subset \E^n\,, F \mbox{ $*$-connected}\,, 
\#F\geq c n^{\frac {(\eps(n)-d)(d-1)}{(\eps(n)-1)d}} \} 
\ge \beta \,.
\end{displaymath} 

Since, for large enough $n$, we have $c n^{\frac {(\eps(n)-d)(d-1)}{(\eps(n)-1)d}}
\geq (\log n)^{\frac 32}$, it might be easier to prove that: 

for $p>p_c$, there exists a constant $\beta$ such that, 
Q.a.s. on the set $\#\CC=\infty$, and for large enough $n$, one has 
\begin{equation}
\label{*3}
\inf \{ 
\frac{\sum_{\{e \in  F\}} {\bf 1}_{\omega(e)=1}}
     {\#F}\,, 
F \subset \E^n\,, F \mbox{ $*$-connected}\,, 
\#F\geq (\log n)^{\frac 32} \} 
\ge \beta \,.
\end{equation}

\subsection{Large values of $p$}
\label{largep}

We use formulation (\ref{*3}), together with a contour argument, 
to check that (\ref{*}) holds when $p$ is close enough 
to $1$. 

For a given set of edges $F\subset\E^n$, since the random variables $(\omega(e),e\in F)$ 
are Bernoulli(p) and independent, we have:
\begin{equation}
\label{3.2}
Q\left[\frac{\sum_{\{e \in  F\}} {\bf 1}_{\omega(e)=1}}
     {\#F}\leq\beta\right]
     \leq
e^{\lambda\beta\#F} (pe^{-\lambda}+(1-p))^{\#F}\,.
\end{equation}

On the other hand, we recall that the number of $*$-connected sets $F\subset\E^n$ of 
cardinality $m$ such that $0\in F$ is bounded by $\exp(am)$, for some constant $a$ that depends on the dimension $d$ only. 
See \cite{sinai82a} for instance. 
Therefore 
\begin{eqnarray*}
&
Q \left[ \exists F\subset\E^n, F \mbox{ $*$-connected}\,, \#F\geq (\log n)^{\frac 32}, 
\mbox{ and}\displaystyle \sum_{\{e \in  F\}} {\bf 1}_{\omega(e)=1}\leq \beta \#F\right]\\  
\leq& \displaystyle \sum_{m\ge (\log n)^{3/2}}
(2n+1)^d \, e^{am+\lambda\beta m} (pe^{-\lambda}+(1-p))^m\,. \\  
\end{eqnarray*}

If $p$ is close enough to $1$, and if we choose $\beta$ small enough, there will exist 
$\lambda >0$ such that: 
\begin{displaymath}
-\xi=a+\lambda\beta+\log((1-p)+e^{-\lambda}p)<0\,.
\end{displaymath}

Then, for large enough $n$, one has: 
\begin{eqnarray*}
Q \big[ \exists F\subset\E^n,  &F \mbox{ $*$-connected}\,,& \#F\geq (\log n)^{\frac 32}, 
\mbox{ such that}  \\ 
&&  \sum_{\{e \in  F\}} {\bf 1}_{\omega(e)=1}\leq \beta \#F\big]
\leq e^{-\frac \xi 2 (\log n)^{3/2}}\,.
\end{eqnarray*}
Note that this last expression is sumable in $n$. 
Therefore the Borel-Cantelli lemma implies that, $Q$.a.s., for large enough $n$, 
for all $*$-connected set $F\subset\E^n$ with $\#F\geq (\log n)^{3/2}$, we have 
$\sum_{\{e \in  F\}} {\bf 1}_{\omega(e)=1}\geq \beta \#F$. In view of 
(\ref{*3}), we deduce that we have proved that there exists a number $p(d)<1$ 
such that (\ref{*}) holds for any $p>p(d)$ and any $\alpha\in]0,\frac 12 ]$. 
\carre

\paragraph{Remark:} for further purposes, we state here estimates of the tail
of the distribution of the Cheeger constant:  

\begin{displaymath}
\II_\infty(\CC^n) 
=
\inf_{A \subset \CC^n\,,\# A \leq \frac {\# \CC^n}2 }
\frac{\# (\partial_{\CC^n} A)}{\# A}\,.
\end{displaymath} 

We wish to estimate $Q[\II_\infty(\CC^n)\leq \frac \beta n]$. 
As in the previous computation, we may restrict our attention to sets $A$ such
that $\#A\geq \frac n \beta$. Also as in previous computation, we see that, on
the event $\#\CC^n\geq 2\alpha (2n+1)^d$, then $\#D\leq
(1-\alpha)(2n+1)^d$. 

Since $\# D\geq\# A\geq n/\beta$,
we have  
\begin{displaymath}
\#F\geq \II^{(\alpha)}(\#D)^{\frac{d-1}d}\geq\II^{(\alpha)}
(\frac n \beta)^{\frac{d-1}d}\,.
\end{displaymath}
On the other hand, 
because $\#A\leq (2n+1)^d$, then 
\begin{displaymath}
\#F\geq \II^{(\alpha)} (\# A)^{\frac{d-1}d}
\geq 
\frac {\II^{(\alpha)}}{2n+1} \#A \,.
\end{displaymath}
(We kept on using the notation $F=\partial_{\BB^n}D$.) 
Thus 
\begin{displaymath}
\frac{\# (\partial_{\CC^n} A)}{\# A}
\geq \frac {\II^{(\alpha)}}{2n+1} 
\frac {\sum_{\{e \in  F\}} {\bf 1}_{\omega(e)=1}}{\#F}\,.
\end{displaymath}

Thus $\#\CC^n\geq 2\alpha (2n+1)^d$ and $\II_\infty(\CC^n)\leq\frac \beta
n$ imply that there exists a $*$-connected set $F\subset\E^n$ such that $\#F\geq
\II^{(\alpha)} (\frac n \beta)^{(d-1)/d}$ and 
\begin{displaymath}
\sum_{\{e \in  F\}} {\bf 1}_{\omega(e)=1}\leq 3\frac  \beta {\II^{(\alpha)}} \#F\,.
 \end{displaymath}

As in (\ref{3.2}), we then get that there exists $p(d)<1$ such that for 
any $\alpha>0$, there exist constants $\beta>0$ and $\xi>0$ such that, for all
$p\in]p(d),1[$, one has: 
\begin{eqnarray*}
Q \left[\II_\infty(\CC^n)\leq \frac \beta n \hbox { and } 
               \#\CC^n\geq 2\alpha (2n+1)^d  \right]
&\leq &\sum_{m\geq \II^{(\alpha)} (n/\beta)^{(d-1)/d}} 
      (2n+1)^d e^{-2\xi m}\\
&\leq & e^{-\xi n^{(d-1)/d}}
\end{eqnarray*}  
for large enough $n$.

We also know, see \cite {deuschel}, 
that  for some constants $\alpha>0$ and $\xi>0$, we have 
\begin{displaymath}
Q\left[\#\CC^n\leq 2\alpha (2n+1)^d  \right] 
\leq e^{-\xi n^{d-1}}\,.
\end{displaymath}

Gathering the last two inequalities, we obtain that:

\begin{theorem}
\label{prop}
For any dimension $d\geq 2$, there exists $p(d)<1$, and constants $\beta>0$
and $\xi>0$ such that, for all $p\in]p(d),1[$ and $n\geq 1$, 
\begin{displaymath}
Q\left[\II_\infty(\CC^n)\leq \frac \beta n\right]
\leq e^{-\xi n^{(d-1)/d}}\,.
\end{displaymath}
\end{theorem}

As a consequence of this estimate and Cheeger's inequality, we obtain a lower
bound for the spectral gap: 
\begin{equation}
Q\left[\lambda^n\leq \frac \beta {n^2}\right]
\leq e^{-\xi n^{(d-1)/d}}\,.
\end{equation}

\subsection{Another isoperimetric inequality}
\label{another}

In order to prepare for the use of renormalization arguments, we shall 
need a more sophisticated version of the inequality (\ref{*}), still for values of $p$ 
close enough to $1$. We shall also have to consider both site and bond percolation 
models. 

As before, we let $\BB^n$ be the box $[-n,n]^d$; $\CC^n$ is the connected component 
of the random graph $\omega$ that contains the origin. We also use the notation 
$\GG^n$ to be the set of vertices of $\BB^n$ that belong to $\omega$. Thus $\CC^n$ is a 
subset of $\GG^n$. Let $\LL^n$ be the largest connected component of $\omega$ in $\BB^n$.

Let $A$ be a subset of $\BB^n$. Define $n(A)$ to be the number of connected components 
of $\BB^n\setminus\LL^n$ that contain at least one  connected component of $A$. 

Let $p<1$, $\alpha\in[\frac 12, 1[$ and $\beta>0$ and define the event $\AAA$ by: 

for large enough $n$, 
\begin{equation}
\label{*4}
\inf \left\{ 
\frac{n(A)+\#(\partial_{\LL^n} A)}
     {(\#A)^{\frac{\eps(n)-1}{\eps(n)}}}\,, 
A \subset \BB^n\,, \#A\leq \alpha\#\BB^n \right\} 
\ge \beta n^{\frac d{\eps(n)}-1} \,.
\end{equation}

We shall prove that there exists $p_0<1$ such that, 
for all $\alpha\in[\frac 12, 1[$,  $p>p_0$ there exists  $\beta>0$ s.t.    
$Q$.a.s.,  $\AAA$ holds {\it i.e.} we prove that 
\begin{equation}
\label{*5} 
Q[\AAA]=1\,.
\end{equation} 

Before entering the proof, let us state some simple remarks.  

We first note that, on the set $\#\CC=\infty$, for large enough $n$, $\LL^n=\CC^n$. 
(On the set $\#\CC=\infty$, then $\#\CC^n\geq 2\alpha(2n+1)^d$ for some $\alpha>0$, 
see Appendix \ref{app-card}. It is easy to show that, for $p$ close enough to $1$, there is 
at most one connected component of $\GG^n$ of size larger than $(\log n)^{2d/(d-1)}$, 
see below.) 
Thus, once we have proved (\ref{*5}),  it will follow that, for some constants $p<1$ and $\beta>0$, $Q$.a.s. on the 
set $\#\CC=\infty$, for large enough $n$, one has: 
\begin{equation}
\label{*6}
\inf \left\{ 
\frac{n(A)+\#(\partial_{\CC^n} A)}
     {(\#A)^{\frac{\eps(n)-1}{\eps(n)}}}\,, 
A \subset \BB^n\,, \#A\leq \alpha\#\BB^n \right\} 
\ge \beta n^{\frac d{\eps(n)}-1} \,.
\end{equation} 
In this last statement, $n(A)$ may as well be defined as 
the number of connected components 
of $\BB^n\setminus\CC^n$ that contain at least one  connected component of $A$. 
If, in (\ref{*6}),  we restrict ourselves to sets $A$ which are contained in $\CC^n$, 
then $n(A)=0$ and 
we retrieve the isoperimetric inequality (\ref{*}). 

$\AAA$ is an increasing event. Indeed, let $A\subset\BB^n$ and 
suppose we add one edge to $\omega$. Then $\#(\partial_{\LL^n} A)$ will not 
decrease. Assume that $n(A)$  decreases by $1$. It implies that at least one 
of the connected components of $A$ did not intersect $\LL^n$  
before the addition of the extra edge and intersects $\LL^n$ after the addition 
of the extra edge. Then a new edge appeared in $\#(\partial_{\LL^n} A)$. We conclude 
that the addition of one edge to $\omega$ does not decrease the sum 
$n(A)+\#(\partial_{\LL^n} A)$. 

As a consequence, we deduce that if (\ref{*5}) holds for some $p_0<1$, it then holds 
for all $p\in[p_0,1]$. From the comparison theorems established in \cite{liggett97}, it is 
equivalent to prove (\ref{*5}) for bond or site percolation.

We now turn to the proof of (\ref{*5}):  
by choosing $p$ close enough to $1$, we may, and will, always assume that 
$\#\LL^n\geq \frac {1+\alpha}2 \#\BB^n$, see Appendix \ref{app-card}. 
Since $\#A\leq\alpha\#\BB^n$, $\LL^n$ cannot be contained in $A$ and 
therefore $n(A)+\#(\partial_{\LL^n} A)\geq 1$. 
Therefore we may, and will, restrict ourselves to sets $A$ such that 
$\#A\ge \beta^{-\frac{\eps(n)}{\eps(n)-1}}n^{\frac {\eps(n)-d}{\eps(n)-1}}$. 
For large enough $n$, 
$\beta^{-\frac{\eps(n)}{\eps(n)-1}}n^{\frac {\eps(n)-d}{\eps(n)-1}}
\geq(\log n)^{{\frac 32}{\frac d{d-1}}}$. 
Therefore 
we may assume that $\#A\ge (\log n)^{ {\frac 3 2} {\frac d{d-1}} }$. 

Let $A_1$ be the union of the connected components of $A$ that intersect $\LL^n$. 
Let $A_2$ be the union of the connected components of $A$ that do not intersect 
$\LL^n$. Note that $\#(\partial_{\LL^n} A)=\#(\partial_{\LL^n} A_1)$. 

Let us first assume that $\#A_1\leq \frac 12 \#A$. Then $\#A_2\geq\frac 12 \#A$ and therefore\\
$\#A_2\geq\frac 12 (\log n)^{{\frac 32}{\frac d{d-1}}}$. 

Let $c_1$, ...$c_k$ be the different connected components of $\BB^n\setminus\LL^n$. 
The same contour argument as in paragraph \ref{largep} shows that, for any $a>1$, 
if we choose $p$ close enough to $1$, then, for large enough $n$, $*$-connected sets 
of volume larger than $(\log n)^a$ intersect $\LL^n$. 
Since, for all $i$, $\partial_{\BB^n} c_i$ is $*$-connected, we have 
$\#\partial_{\BB^n} c_i\leq (\log n)^a$. The classical isoperimetric inequality 
then implies that  
$\#c_i\leq (\log n)^{\frac{ad}{d-1}}$, for all $i=1...k$. 
Here we choose $a=2$.  

We then have:
\begin{eqnarray*}
n(A)&=&\sum_{i=1}^k 1_{A_2\cap c_i\not=\emptyset}\\
&\geq&\sum_{i=1}^k \frac{\#(A_2\cap c_i)}{(\log n)^{\frac{2d}{d-1}}}\\
&=& \frac{\#A_2}{(\log n)^{\frac{2d}{d-1}}}\\
&\geq& n^{\frac d{\eps(n)}-1} (\# A_2)^{\frac{\eps(n)-1}{\eps(n)}}\,,
\end{eqnarray*}
where the last inequality comes from the fact that 
$n^{1-\frac d{\eps(n)}}\geq {(\log n)^{\frac{2d}{d-1}}}$ for large $n$. 
We conclude that $n(A)\ge  (\frac{\# A}2)^{\frac {\eps(n)-1}{\eps(n)}} n^{\frac d{\eps(n)}-1}$.  

Now assume that $\#A_1\geq \frac 12 \#A$. Then  
$\#A_1\geq\frac 12 (\log n)^{{\frac 32}{\frac d{d-1}}}$. 
We wish to show that 
\begin{equation}
\label{*7}
\frac{\#(\partial_{\LL^n} A_1)}
     {(\#A_1)^{\frac{\eps(n)-1}{\eps(n)}}}\geq \beta n^{\frac d{\eps(n)}-1}\,.
\end{equation}

Without loss of generality, we may assume that $A_1$ is connected. 
Let us first check that we may also assume that all the connected components 
of $\BB^n\setminus A_1$ intersect $\LL^n$ and have size bounded by $\frac 12 \#\BB^n$. 
Indeed, let $c_1, ..., c_k$ be the connected components of $\BB^n\setminus A_1$. 
Suppose that some $c_i$ does not intersect $\LL^n$.   
Then $\#A_1\leq\#(A_1\cup c_i)$  and $\#\partial_{\LL^n}A_1=\partial_{\LL^n}(A_1\cup c_i)$. 
Also note that $\#(A_1\cup c_i)\le \frac{1+\alpha}2\#\BB^n$. 
Therefore, changing the value of $\alpha$ to $\frac{1+\alpha}2$,  it is sufficient 
to prove (\ref{*7}) for $A_1\cup c_i$ instead of $A_1$. 
From now on, we shall assume that all $c_i$ intersect $\LL^n$. 

Assume that $\#c_i>\frac 12 \#\BB^n$ for some $i$. Let $A_3=\BB^n\setminus c_i$. Note that 
$A_3$ is connected, $\BB^n\setminus A_3$ is connected, 
$A_3$ intersects $\LL^n$ (because $A_1$ intersects $\LL^n$) 
and $\#A_3\le \frac 12 \#\BB^n$. Since $\partial_{\LL^n} A_3\subset\partial_{\LL^n} A_1$ and 
$A_1\subset A_3$, if we can prove (\ref{*7}) for $A_3$, we can prove it for $A_1$. 
Thus it is no loss of generality to assume that  $\#c_i\le \frac 12 \#\BB^n$. 

We shall now use the fact that, for some choice of $p$ and $\beta$, we have 
\begin{equation}
\label{*8}
\frac{\#\partial_{\LL^n}c_i}{(\#c_i)^{\frac{\eps(n)-1}{\eps(n)}}}
\geq\beta n^{\frac d{\eps(n)}-1}\,.
\end{equation}
Equation (\ref{*8}) follows from the fact that  $c_i$ is connected, 
$\BB^n\setminus c_i$ is connected and $\partial_{\LL^n}c_i\not=\emptyset$, 
since $c_i$ intersects $\LL^n$. Therefore the results 
of paragraph \ref{largep} can be applied. 

We conclude by noticing that: 

\begin{eqnarray*}
\#\partial_{\LL^n}A_1 &=& \sum_i \#\partial_{\LL^n} c_i\\
&\geq& \beta n^{\frac d{\eps(n)}-1} \sum_i (\# c_i)^{\frac{\eps(n)-1}{\eps(n)}}\\
&\geq& \beta n^{\frac d{\eps(n)}-1} (\sum_i \# c_i)^{\frac{\eps(n)-1}{\eps(n)}}\\
&=& \beta n^{\frac d{\eps(n)}-1} (\#\BB^n-\#A_1)^{\frac{\eps(n)-1}{\eps(n)}}\\
&\geq& \beta n^{\frac d{\eps(n)}-1} (\frac{1-\alpha}\alpha)^{\frac{\eps(n)-1}{\eps(n)}} 
(\#A_1)^{\frac{\eps(n)-1}{\eps(n)}}\,,
\end{eqnarray*}
since $\#A\le\alpha\#\BB^n$. 
\carre

\subsection{Renormalization}
\label{renormalization}

We now explain how to push the isoperimetric inequality from large values of $p$ 
down to any $p>p_c$. To this end, we mainly rely on Proposition 2.1. in 
\cite{antal}. We shall also use some terminology from \cite{antal}. 

Choose $p>p_c$. $N$ is an integer. We chop $\Z^d$ in a disjoint union of boxes of side 
length $2N+1$. Say $\Z^d=\cup_{{\mathbf i}\in\Z^d}B_{\mathbf i}$, where $B_{\mathbf i}$ is the box of 
center $(2N+1){\mathbf i}$. Still following \cite{antal}, let $B'_{\mathbf i}$
be the box of center $(2N+1){\mathbf i}$ and side length $\frac 5 2 N +1$. 
 From now on, the word {\it box} will mean one of the boxes 
$B_{\mathbf i}, {\mathbf i}\in\Z^d$. 

We say that a box $B$ is {\it good} if 
$B$ contains at least one edge from $\omega$ 
and the event $R_{\mathbf i}^{(N)}$ in equation (2.9) of 
\cite{antal} is satisfied. Otherwise, $B$ is a {\it bad} box. 
We recall that the event $R_{\mathbf i}^{(N)}$ is defined by: 
there is a unique cluster, $K$, in $B'_{\mathbf i}$; all open paths 
contained in $B'_{\mathbf i}$ and of radius larger than $\frac 1 {10} N$ 
intersect $K$ within $B'_{\mathbf i}$; $K$ is crossing for each subbox 
$B\subset B'_{\mathbf i}$ of side larger than $\frac 1 {10} N$. 
See \cite{antal} for a more precise definition. 

Let $A$ be a subset of $\CC^n$. We say that $A$ {\it touches} the box $B$ if 
$A\cap B\not=\emptyset$. We say that $A$ {\it fills} $B$ if $A$ touches $B$ and 
$A\cap B=B\cap\CC^n$. 
We use the notation $\bar n_1$ (resp $n_1$) to denote the number of boxes 
(resp. good boxes) touched by $A$ but not filled by $A$. Similarly, 
let $\bar n_2$ (resp $n_2$) denote the number of boxes 
(resp. good boxes) filled by $A$. 

Following \cite{antal}, we call {\it renormalized} process the percolation model obtained 
by taking the image of the initial percolation model by the application $\phi_N$, 
see equation (2.11) in \cite{antal}. A site $\mathbf i \in\Z^d$ is thus declared 
{\it white} if the box $B_{\mathbf i}$ is good. 

We choose $N$ large enough so that the following two conditions are satisfied: 
\begin{enumerate}
\item the isoperimetric inequality (\ref{*4}) holds for the renormalized percolation model and 
$\alpha=\frac 4 7 (1+10^{-3})$, 
\item for large enough $n$ and for all connected set $A\subset \BB^n$ with \\
$\# A\geq (2N+1)^d(\log n)^{3/2}$, we have 
\begin{displaymath}
n_1+n_2\geq \frac 7 8 ({\bar n}_1+{\bar n}_2) \,.
\end{displaymath}
\end{enumerate}

Such $N$ exists. Indeed, we already know from Section \ref{another} than $p'$ can be chosen 
close enough to $1$ so that the isoperimetric inequality (\ref{*4}) is satisfied for site 
percolation of parameter $p'$. From the comparison of Proposition 2.1 in \cite{antal}, 
and the fact that event $\AAA$ in (\ref{*4}) is increasing, we 
deduce that (\ref{*5}) is also satisfied  for the 
renormalized process if $N$ is large enough. 

Let us check that we can choose $N$ so that {\rm({\it ii})} is satisfied. For a connected
set $A\subset\BB^n$ with $\# A\geq (2N+1)^d(\log n)^{3/2}$, we have 
\begin{displaymath}
{\bar n}_1+{\bar n}_2\geq (2N+1)^{-d}\#A\geq (\log n)^{3/2}\,. 
\end{displaymath}
Call $\tilde A$,  the set of indices $\mathbf i$ such that $B_{\mathbf i}$ is touched by $A$. 
Note that ${\bar n}_1+{\bar n}_2=\#\tilde A$ and that $\tilde A$ is connected. 
Therefore, using Proposition 2.1 in \cite{antal}, we only have
to check the following property: 

 for any constant $c$, for $p$ close enough to $1$, and for site percolation,
 $Q$.a.s., for large enough $n$, for any connected set ${\tilde A}\subset \BB^n$ with
 $\#{\tilde A}\geq c(\log n)^{3/2}$, then 
\begin{equation}
\label{fuck}
\#({\tilde A}\cap\GG^n)\geq \frac 7 8 \#{\tilde A}\,.
\end{equation}
(Remember than $\GG^n$ is the set of vertices of $\omega$ in the box $\BB^n$.) 

Inequality (\ref{fuck}) follows from a Borel-Cantelli 
argument based on the fact that the number of connected subsets of $\Z^d$ containing $0$ and 
of cardinal $m$ is bounded by $\exp(a m)$, 
for some constant $a$ that depends on the dimension only.  
\carre

We want to check that claim (\ref{*2}) holds true.

First let  $A$ be a connected subset of $\CC^n$ such that 
$\# A\geq c n^{(\eps(n)-d)/(\eps(n)-1)}$ and 
$n_1\geq 2^d\beta (\#A)^{(\eps(n)-1)/\eps(n)} n^{d/\eps(n)-1}$. \\
Note that $cn^{(\eps(n)-d)/(\eps(n)-1)}\geq (2N+1)^d(\log n)^{3/2}$ if $n$ is large
enough. It follows that $A$ is not entirely contained in one single box. 

Each good box touched, but not filled, by $A$ contributes by at least one edge 
the boundary of $A$. Indeed, assume that $A$ touches but does not fill
$B_{\mathbf i}$. Then we can find two points $x,y\in B_{\mathbf i}\cap\CC^n$
such that $x\in A$ and $y\notin A$. Since $A$ is not entirely contained in
$B_{\mathbf i}$, there is
a path in $\omega$ of length at least $N/5$ linking $x$ to the boundary of
$B_{\mathbf i}$. Similarly, there is a path in $\omega$ of length at least $N/5$
linking $y$ to the boundary of $B_{\mathbf i}$. From the definition of a good box, 
it follows that these two open paths have to
be connected to each other within $B'_{\mathbf i}$. Therefore, there is an
open path linking $x$ to $y$ within $B'_{\mathbf i}$. On this open path, there
is an edge, say $(a,b)$, where $a$ and $b$ are neighbors, both are in $\CC^n$
and $a\in A$, $b\notin A$. Thus we have found one edge in $\partial
_{\CC^n}A$. We can repeat that contruction for each of the $n_1$ good boxes
touched but not filled by $A$. Since a given edge cannot appear more that
$2^d$ times, we get that  $\#(\partial_{\CC^n}A)\geq 2^{-d} n_1$. 
Therefore, since $n_1\geq 2^d\beta (\#A)^{(\eps(n)-1)/\eps(n)} n^{d/\eps(n)-1}$, then 
we also have 
$\#(\partial_{\CC^n}A)\geq \beta (\#A)^{(\eps(n)-1)/\eps(n)} n^{d/\eps(n)-1}$.

We let now $A$ be a connected subset of $\CC^n$  with 
$n_2\leq \frac 12(\frac{2n+1}{2N+1})^d$, \par 
$n_1\leq 2^d\beta (\#A)^{(\eps(n)-1)/\eps(n)} n^{d/\eps(n)-1}$ 
and $\# A\geq c n^{(\eps(n)-d)/(\eps(n)-1)}$. Note that $c
n^{(\eps(n)-d)/(\eps(n)-1)}\geq (2N+1)^d(\log n)^{3/2}$ if $n$ is large
enough. It follows that $A$ is not entirely contained in one single box.

Let ${\tilde A}$ be the set of indices ${\mathbf i}\in\Z^d$ such that $B_{\mathbf i}$ is touched 
by $A$. We wish to use the isoperimetric inequality (\ref{*4}) for the set ${\tilde A}$. 
By definition $({\bar n}_1+{\bar n}_2)=\# {\tilde A}$. 

Since $({\bar n}_1+{\bar n}_2)(2N+1)^d\geq \#A$ and 
$\frac 7 8({\bar n}_1+{\bar n}_2)\leq n_1+n_2$, we then have 
\begin{displaymath}
n_2\geq (2N+1)^{-d} \frac 7 8 \# A-2^d\beta n^{\frac d{\eps(n)}-1}
    (\# A)^{\frac{\eps(n)-1}{\eps(n)}}
\geq (2N+1)^{-d} \frac 6 8 \# A \,,
\end{displaymath}
for large enough $n$. In particular it follows that, for large enough $n$, 
$n_1\le 10^{-3} n_2$. 

Remember that we have assumed that 
$n_2\leq \frac 1 2 (\frac {2n+1}{2N+1})^d $. 
Thus 
\begin{displaymath}
({\bar n}_1+{\bar n}_2)
\leq \frac 4 7 (1+10^{-3})(\frac {2n+1}{2N+1})^d\,.
\end{displaymath}

Therefore we may apply the isoperimetric inequality (\ref{*4}) to the set $\tilde A$ 
with $\alpha=\frac 4 7 (1+10^{-3})$. 
Clearly $\tilde A$ is connected. We use the notation $\tilde .$ to indicate quantities 
defined at the level of the renormalized process. 

Either $\tilde A$ is contained in one of the connected components 
of ${\tilde \BB^n}\setminus {\tilde \LL^n}$. Then ${\tilde n}({\tilde A})=1$ and 
$\partial_{\tilde \LL^n}{\tilde A}=\emptyset$ and therefore 
$1\geq \beta n^{d/\eps(n)-1} (\#{\tilde A})^{(\eps(n)-1)/\eps(n)}$. Since 
$\partial_{\CC^n}A\not=\emptyset$, we then have 
\begin{displaymath} 
\# (\partial_{\CC^n} A)\geq 
\beta n^{d/\eps(n)-1} (\#{\tilde A})^{(\eps(n)-1)/\eps(n)}\,.
\end{displaymath}

Let us now suppose that $\tilde A$ intersects $\tilde \LL^n$. Then 
${\tilde n}({\tilde A})=0$ and 
\begin{displaymath}
\#(\partial_{\tilde\LL^n}{\tilde A})\geq 
\beta n^{\frac d{\eps(n)}-1} (\#{\tilde A})^{\frac{\eps(n)-1}{\eps(n)}}\,.
\end{displaymath}

Let $({\mathbf i},{\mathbf i'})\in \partial_{\tilde\LL^n}{\tilde A}$. Then  
$B_{\mathbf i}$ and $B_{\mathbf i'}$ are good boxes; 
${\mathbf i}$ and ${\mathbf i'}$ are neighbors; 
$\tilde A$ intersects $B_{\mathbf i}$ but not $B_{\mathbf i'}$.  
But note that each such couple $({\mathbf i},{\mathbf i'})$ contributes by at least 
one edge to the boundary 
of $A$ in $\CC^n$. Indeed, we can find points $x\in B_{\mathbf i}$ and $y\in
B_{\mathbf i'}$ such that $x\in A$ and $y\in\CC^n\setminus A$. Because the two
boxes $B_{\mathbf i}$ and $B_{\mathbf i'}$ are good, there is an open path
linking $x$ to $y$ within $B'_{\mathbf i}\cup B'_{\mathbf i'}$, and, on this
path, there must be an edge of $\partial_{\CC^n}A$. As we perform this
construction for different choices of $({\mathbf i},{\mathbf i'})$, a 
given edge appears at most $2^d$ times. 
Therefore 
\begin{displaymath}
\#(\partial_{\CC^n} A)
\geq 2^{-d}\beta  n^{\frac d{\eps(n)}-1} (\#{\tilde A})^{\frac{\eps(n)-1}{\eps(n)}}\,.
\end{displaymath} 

Since $\#{\tilde A}\geq n_2\geq (2N+1)^{-d} \frac 6 8 \# A$, we have  
\begin{displaymath} 
\#(\partial_{\CC^n} A)
\geq \beta n^{\frac d{\eps(n)}-1} (\# A) ^{\frac{\eps(n)-1}{\eps(n)}}\,,
\end{displaymath} 
with a different value of $\beta$. 

Finally, let $A$ be a connected subset of $\CC^n$ such that 
$\ A'=\CC^n\setminus A$ is also connected, $\#A\leq\#\CC^n/2$, 
$\ n_2 \geq \frac 12(\frac{2n+1}{2N+1})^d$, 
$\ n_1 \leq 2^d\,\beta\, (\#A)^{(\eps(n)-1)/\eps(n)}\, n^{d/\eps(n)-1}$ 
and \\$\ \# A\geq c n^{(\eps(n)-d)/(\eps(n)-1)}$. 
Let $n_1'$ and $n'_2$ be defined as $n_1$ and $n_2$ with $A$ being replaced 
by $A'$. Note that $n_1'=n_1\leq 2^d\beta (\#A)^{(\eps(n)-1)/\eps(n)} n^{d/\eps(n)-1}$.\\
Also $\#A'\geq \#\CC^n/2\geq c n^{(\eps(n)-d)/(\eps(n)-1)}$, see Appendix \ref{app-card}. 
Since $n_2+n'_2\leq (\frac{2n+1}{2N+1})^d$, we must have 
$n_2'\leq \frac 12(\frac{2n+1}{2N+1})^d$. 
Thus we may apply the previous isoperimetric inequality to $A'$ and get that: 
\begin{displaymath} 
\#(\partial_{\CC^n} A')
\geq \beta n^{\frac d{\eps(n)}-1} (\# A') ^{\frac{\eps(n)-1}{\eps(n)}}\,.
\end{displaymath} 
But note that $\partial_{\CC^n} A'=\partial_{\CC^n} A$ and 
$\#A'\geq\frac 12 \#\CC^n\geq \#A$. Therefore: 
\begin{displaymath} 
\#(\partial_{\CC^n} A)
\geq \beta n^{\frac d{\eps(n)}-1} (\# A) ^{\frac{\eps(n)-1}{\eps(n)}}\,.
\end{displaymath}

\carre

\section{Proofs of the theorems}

We prove inequalities (\ref{ineq-th1}) and 
(\ref{ineq-th2}),{\it i.e.} that 
\begin{displaymath}
 \sup_{x \in \CC} P_0^{\omega}[X_t=x] 
 \leq 
 \frac{c_1}{t^{d/2}}\,.
\end{displaymath}
We use Carne-Varopoulos inequality (cf. inequality (\ref{carne}), Appendix \ref{carn-var}):
\begin{displaymath}
 P_0^{\omega}[X_t=x]
 \leq
 e^{-\frac{|x|^2}{4\,t}}+e^{-ct}\,.
\end{displaymath}
\paragraph{First case:}  If $|x|^2 \geq 2 d\,t\,\log t$, then 
$e^{-\frac{|x|^2}{4\,t}}\leq C t^{-d/2}$ and also $e^{-ct}\leq C t^{-d/2}$,
provided $t$ is large enough.

\paragraph{Second case:} If $|x|^2 < 2d\,t\,\log t$.

 We choose $t \log t=b\,n^2$, with $b < \frac{1}{4d+2}$.

 Let $\tau^n$ denote the exit time of $\BB^{n-1}$ for $X_t$. Then: 
\begin{displaymath}
 P_0^{\omega}[X_t=x] 
 \leq 
 P_0^{\omega}[X^n_t=x]+P_0^{\omega}[\tau^n \leq t]\,.
\end{displaymath}

Combining with the estimates (\ref{estim-reflected}) and  (\ref{exit}), we get that: 
\begin{displaymath}
 P_0^{\omega}[X_t=x] 
 \leq 
 n^{-d} + 
\left(\frac{4\,\eps}{\beta^2}\right)^{\frac{\eps}{2}}
\,\frac{n^{\eps-d}}{t^\frac{\eps}{2}}
 +2 t\,n^{d-1} e^{-\frac{n^2}{4\,t}}+e^{-ct}\,.
\end{displaymath}
We consider the right side of the inequality term by term:
\begin{enumerate}
\item 
 $n^{-d}=\left( \frac{t\,\log t}{b}\right)^{-d/2} \leq t^{-d/2},$ this 
will be satisfied as soon as $t$ is large enough.
\item
\begin{eqnarray*}
\left(\frac{4\,\eps}{\beta^2}\right)^{\frac{\eps}{2}}
\,\frac{n^{\eps-d}}{t^\frac{\eps}{2}}
\leq \frac{C}{t^\frac{d}{2}}
&\Leftrightarrow&
\left(\frac{4\,\eps}{\beta^2}\right)^{\frac{\eps}{2}}
\,\left(\frac{n^2}{t}\right)^{\frac{\eps-d}{2}}
\leq C
\\
&\Leftrightarrow&
\left(\frac{4\,\eps}{\beta^2}\right)^{\frac{\eps}{2}}
\,\left(\frac{\log t}{b}\right)^{\frac{\eps-d}{2}}
\leq C
\\
&\Leftrightarrow&
\log K+d\,\frac{\log \log n}{\log n}\,\log \left(\frac{\log t}{b}\right)
\leq \log C\,,
\end{eqnarray*}
with $K=\left(\frac{4\,\eps}{\beta^2}\right)^{\frac{\eps}{2}}$.
Since $\frac{\log \log n}{\log n}\,\log\log t \rightarrow 0$, this is true for
large $t$, and for some constant $C$.
\item
\begin{eqnarray*}
 2t\,n^{d-1} e^{-\frac{n^2}{4\,t}}
 \leq \frac{C}{t^\frac{d}{2}}
&\Leftrightarrow&
\log t+ (d-1)\,\log n+ \frac{d}{2}\,\log t \leq \log \frac C 2 + \frac{n^2}{4\,t}
 \\
&\Leftrightarrow&
\log t+ (d-1)\,\log n+ \frac{d}{2}\,\log t \leq \log \frac C 2 + \frac{\log
t}{4\,b}
\\
&\Leftrightarrow&
\log t \left[ \frac{1}{4b}-1-\frac{d}{2}\right]
+ \log \frac C 2 \geq (d-1)\,\log n\,,
\end{eqnarray*}
and this is OK if $2\,\left[ \frac{1}{4b}-1-\frac{d}{2}\right] >
d-1$,{\it i.e.} if $b$ satisfies $b < \frac{1}{4d+2}$\,.
\item
The last term, $e^{-ct}$, clearly decays faster than $t^{-d/2}$.
\end{enumerate}

\carre

\newpage
\appendix
\section*{Appendix}

\section{Connectivity}
\label{app-con}

We prove the following statement:

\begin{quote}
 Let $\BB$ be a finite box of $\Z^d$. Let $A \subset \BB$ such that $A$ and
 $\BB\setminus A$ are both connected.
 Then, $\partial_{\BB} A$ is  $*$-connected.
\end{quote} 

Let us recall what $*$-connectedness is,  
see \cite{deuschel}: we embed $\Z^d$ into $\R^d$ and endow 
$\R^d$ with the norm $\Vert x\Vert=\max_{i=1...d}\vert x_i\vert$. Consider two nearest 
neighbor edges in $\Z^d$, say $e=(x,y)$ and $e'=(x',y')$, where $x\sim y$ and $x'\sim y'$. 
We say that $e$ and $e'$ are $*$-neighbors if 
$\Vert \frac{x+y}2-\frac{x'+y'}2\Vert\leq 1$. A set of edges, say $F$, is said to be 
$*$-connected if for any $e\in F$ and $e'\in F$, there exists a sequence of edges 
starting at $e$, ending at $e'$ and such that any two successive edges in this sequence are 
$*$-neighbors.

\proof

Remember that 
$\partial_{\BB} A = 
  \{ (x,y)\,;\,x \sim y,\ x \in A \mbox{ and } y \in \BB\setminus A \mbox{ or } 
      x \in \BB \setminus A \mbox{ and } y \in A \}$.
In particular, we do not take into account edges with an endpoint outside $\BB$.

 We identify a point $x \in \Z^d$ with the unit cube of $\R^d$ whose center is $x$.
Similarly, we identify an edge $e=(x,y)$ with the common face of the cubes defined by $x$
and $y$.
By this way, $\BB$ is identified with a cube in $\R^d$, say $\tilde{\BB}$.
$A$ is identify with a connected subset of $\tilde{\BB}$, say $\tilde{A}$.
$\tilde{\BB} \setminus \tilde{A}$ is connected, and $\partial_{\BB} A$ is identified with
the boundary of $\tilde{A}$ in $\tilde{\BB}$, say $\partial_{\tilde{\BB}} \tilde{A}$.

We shall prove that $\partial_{\tilde{\BB}} \tilde{A}$ is connected.

   Let $e=(a,b)$ and $ e'=(a',b')$ be edges in $\partial_{\BB} A$. Assume
that $a \in A$, $a' \in A$, $b \in \BB \setminus A$, $b' \in \BB \setminus A$. Since $A$ is
connected, there is a path in $A$ linking $a$ with $a'$. In the same way, since 
$\BB \setminus A$ is connected, there is a path in $\BB \setminus A$ linking $b$ with $b'$.
Joining these two paths and the segments $[a,a']$ and $[b,b']$, 
we obtain a closed loop in $\tilde{\BB}$, say $\gamma$, that contains points
both in the interior of $\tilde{A}$ and $\tilde{\BB} \setminus \tilde{A}$, and with exactly
two points of intersection with $\partial_{\tilde{\BB}} \tilde{A}$, say $\alpha$ and 
$\alpha'$. 
    
    Let $(\gamma(t), t \in [0,1])$ be a smooth deformation of $\gamma$ such that
$\gamma(0)=\gamma$, $\gamma(1)$ is reduced to one point and lies in the interior of
$\tilde{A}$. Further, assume that $\gamma(t)\subset \tilde{\BB}$ for all $t \in [0,1]$.

 Consider the evolution of the set $\gamma(t) \cap \partial_{\tilde{\BB}} \tilde{A}$, and
 note that, as $t$ goes from $0$ to $1$, $\gamma(t) \cap \partial_{\tilde{\BB}} \tilde{A}$
 is made of a finite number of points which are continuously moving on 
 $\partial_{\tilde{\BB}} \tilde{A}$.
 Two points may meet, and then disappear. A new point may also appear, and then split into
 two new points that start wandering around.

Call $(\alpha_t, t \in [0,1])$ the trajectory in $\gamma(t) \cap 
\partial_{\tilde{\BB}} \tilde{A}$ issued from $\alpha$. We thus see that $\alpha_t$ cannot
disappear from $\gamma(t) \cap \partial_{\tilde{\BB}} \tilde{A}$ unless it meets another
point in $\gamma(t) \cap \partial_{\tilde{\BB}} \tilde{A}$. 
A moment of thought should convince the reader that this fact implies that there is an arc
in $\cup_{t \in [0,1]} \gamma(t) \cap \partial_{\tilde{\BB}} \tilde{A}$ joining $\alpha$ and
$\alpha '$.

Thus, $\partial_{\tilde{\BB}} \tilde{A}$ is connected, which implies that 
$\partial_{\BB} A$ is $*$-connected, see \cite{deuschel}, Appendix A. 

\carre


\section{Cardinal of the cluster in the box}
\label{app-card}

We establish that, for all $p>p_c$:
$\exists \alpha$ such that $Q$ a.s. on the set  $\# \CC=+\infty$, one has, for large 
enough $n$:
\begin{equation}
\label{card-cluster}
 \# \CC^n \geq 2\alpha \, (2\,n+1)^d\,.
\end{equation}

\medskip

Choose $\rho\geq 1$. For $x\in\CC$, let $D(0,x)$ denote the 
minimal length of an open path in $\CC$ connecting $0$ and $x$. 
Assume that there exists $x\in\CC\setminus\CC^n$ such that $\vert x\vert\le
n/\rho$. Without loss of generality, we may, and will, assume that $\vert
x\vert\geq n/(2\rho)$. 
Thus the shortest path in $\CC$ linking $0$ to $x$ must leave the box
$\BB^n$ and therefore $D(0,x)\geq n\geq \rho \vert x\vert$. 
From Antal and Pisztora \cite{antal}, theorem 1.1, we know that there exist a
choice of $\rho$ and a 
constant $\beta>0$ such that 
\begin{displaymath}
Q\left[ x\in\CC,\ D(0,x)\geq \rho\vert x\vert\right] \leq e^{-\beta \vert
x\vert}\,,
\end{displaymath}
as $\vert x\vert\rightarrow +\infty$. In particular, 
\begin{eqnarray*}
Q\left[ \exists x\in\CC\setminus\CC^n,\ \vert x\vert\leq\frac n\rho \right]
&\leq&Q\left[ \exists x\in\CC,\ \frac n{2\rho}\leq\vert x\vert\leq \frac n
\rho,\ D(0,x)\geq \rho\vert x\vert\right]\\
&\leq& \sum_{\vert x\vert \in[\frac n{2\rho},\frac n\rho]} 
 e^{-\beta\vert x\vert}\\
&\leq& (2n+1)^d e^{-\frac{\beta n}{2\rho}}\,.
\end{eqnarray*}

From Borel-Cantelli Lemma, we deduce that, Q.a.s, for large enough $n$, 
\begin{displaymath}
\CC\cap\BB^{n/ \rho}\subset\CC^n\,.
\end{displaymath}

It directly follows from the ergodic theorem that 
$\#(\CC\cap\BB^{ n/\rho})/(2n+1)^d$ has an almost sure non
vanishing limit. Therefore $\#\CC^n/(2n+1)^d$ has an almost sure non vanishing
lim inf. 

\carre


\section{Carne-Varopoulos bound}
\label{carn-var}

For a given sub-graph of $\Z^d$, say $\omega$, $P_0^\omega$ denotes the law of
the continuous time random walk on $\omega$ started at $0$,{\it i.e.}, under $P_0^\omega$, the
coordinate process $(X_t,t\geq 0)$ waits for an exponential time of parameter
$1$, then chooses uniformly at random one of its neighbors, say $y$, and moves
to $y$ if $y\in\omega$. Otherwise, $X$ stays still. 

We can also construct $X_t$ as a time changed discrete parameter random walk
on $\omega$. Then $X_t=Y_{N_t}$, where $N_t$ is a Poisson process of parameter
$1$ and $(Y_k,k\in\N)$ is the discrete time ramdom walk on $\omega$ defined by
successively choosing, uniformly at random, one neighbor of the current position and
moving to it if it belongs to $\omega$. 

From \cite{carne85a}, we know that 
\begin{displaymath}
P^\omega_0[Y_k=x]\leq e^{-\frac{\vert x\vert^2}{2k}}
\end{displaymath}

Therefore, 
\begin{eqnarray}
\label{carne}
P^\omega_0[X_t=x]&\leq&  e^{-\frac{\vert x\vert^2}{4t}}+P^\omega_0[N_t\geq 2t]\\
\nonumber
&\leq& e^{-\frac{\vert x\vert^2}{4t}}+e^{-ct}\,,
\end{eqnarray}
where $c=\log 4 -1$. 

Let now $\tau^n$ be the exit time of $X$ from $\BB^{n-1}$. Thus
$\sigma^n=N_{\tau^n}$, where $\sigma^n$ is the exit time for the process $Y$. 
Then: 
\begin{eqnarray*}
 P_0^{\omega}[\sigma^n \leq k]
 &=&
 \sum_{i=0}^k P_0^{\omega}[\sigma^n = i]
\leq
 \sum_{i=0}^k \sum_{y \in  \BB^n \setminus \BB^{n-1}} P_0^{\omega}[Y_i=y]
\\
 &\leq& 
 \sum_{i=0}^k  \sum_{y \in  \BB^n \setminus \BB^{n-1}} e^{-\frac{|y|^2}{2\,i}}
\end{eqnarray*}
(with Carne-Varopoulos' inequality).
Now, as $y \in  \BB^n\setminus \BB^{n-1}\ \Rightarrow \ |y|=n$, we obtain the 
following upper bound:
\begin{displaymath}
 P_0^{\omega}[\sigma^n \leq k]
 \leq
 k\,n^{d-1} e^{-\frac{n^2}{2\,k}}\,.
\end{displaymath}

Thus: 
\begin{eqnarray}
\nonumber
P_0^\omega[\tau^n\leq t] &\leq& 
P_0^\omega[\sigma^n\leq 2t]+P_0^\omega[N_t\geq 2t]\\ 
\label{exit}
&\leq& 2t n^{d-1}e^{-\frac{n^2}{4t}}+e^{-ct}\,,
\end{eqnarray}
for some constant $c>0$. 

\carre 


\section{Lower bound for $P^\omega_0[X_t=0]$.}
\label{app-lower}

To conclude, we briefly discuss the lower bound issue. 
We shall only consider the mean: $Q[P^\omega_0[X_t=0]\vert\#\CC=\infty]$.  
$Q$ a.s. lower bounds{\it i.e.} almost sure lower bounds for 
$P^\omega_0[X_t=0]$ are being investigated at the present time. 
Our result is: 

for all $p>p_c$, there exists a constant $c$ such that, for all $t>0$, 
\begin{equation}
\label{low}
Q[P^\omega_0[X_t=0]\vert\#\CC=\infty]\geq\frac c{t^{d/2}}\,.
\end{equation}

For $p>p_c$, there is, with $Q$ probability one, a unique infinite cluster in $\omega$, 
say $\GG$ and 
$Q[P^\omega_0[X_t=0]\vert\#\CC=\infty]=Q[P^\omega_0[X_t=0]\vert 0\in\GG]$. 

In the next sequence of inequalities, we use the reversibility of the process 
$X_t$ and the translation invariance, {\it i.e.} the fact that 
$Q[P^\omega_x[X_t=x]; x\in\GG]$ does not depend on $x$. 

\begin{eqnarray*}
Q[P^\omega_0[X_t=y]; 0\in\GG]
&=& Q[P^\omega_0[X_t=y]; 0,y\in\GG]\\
&=& Q[\sum_z P^\omega_0[X_{t/2}=z] P^\omega_z[X_{t/2}=y]; 0,y\in\GG]\\ 
&=& Q[\sum_z P^\omega_0[X_{t/2}=z] P^\omega_y[X_{t/2}=z]; 0,y\in\GG]\\ 
&\leq& Q[\sqrt {  \sum_z P^\omega_0[X_{t/2}=z]^2 }
         \sqrt{  \sum_z P^\omega_y[X_{t/2}=z]^2 }; 0,y\in\GG]\\ 
&=& Q[\sqrt{P^\omega_0[X_t=0]}\sqrt{P^\omega_y[X_t=y]}; 0,y\in\GG]\\
&\leq& \sqrt{Q[P^\omega_0[X_t=0]; 0\in\GG]} \sqrt{ Q[P^\omega_y[X_t=y]; y\in\GG]}\\
&=& Q[P^\omega_0[X_t=0]; 0\in\GG]\,.
\end{eqnarray*} 
           
From the invariance principle, see \cite{demasiferrari}, it follows that 
there exists $a>0$ such that, for all $t>0$, we have 
\begin{displaymath}
Q[P^\omega_0[X_t\in \BB^{a\sqrt{t}}]\vert 0\in\GG]\geq \frac 12\,.
\end{displaymath}
In particular, 
\begin{displaymath}
\sup_{y\in\BB^{a\sqrt{t}} }Q[P^\omega_0[X_t=y]\vert 0\in\GG]\geq c t^{-d/2}\,,
\end{displaymath} 
for some constant $c$ that depends only on the dimension.\\ 
Since $\sup_{y\in\BB^{a\sqrt{t}} }Q[P^\omega_0[X_t=y]\vert 0\in\GG]=
        Q[P^\omega_0[X_t=0]\vert 0\in\GG]$, 
(\ref{low}) is proved.
\carre

\bigskip

{\it Acknowledgments:} this work was partially supported by 
FAPESP grants no.99/11109-4 for E.R. and  no.99/0961-1 for P.M. 
We thank the people at IME, S\~ao Paulo for their kind hospitality 
and Enrique Andjel for useful references.

\newpage

\end{document}